\documentclass[english]{article}
\usepackage[T1]{fontenc}
\usepackage[utf8]{inputenc}
\usepackage{morefloats}
\usepackage[letterpaper]{geometry}
\usepackage{float}
\usepackage{amsmath}
\usepackage{tikz}

\usepackage{graphicx}
\usepackage{color}
\usepackage{adjustbox}
\usepackage{diagbox}
\usepackage{multirow}
\usepackage{array}
\usepackage{mathrsfs}
\usepackage{amsmath, amssymb, amsthm}
\usepackage{subfigure}
\usepackage{subfiles}
\usepackage{tikz}

\usepackage{epstopdf}
\graphicspath{{./figures/}}

\newtheorem{theorem}{Theorem} 

\newtheorem{assumption}[theorem]{Assumption}
\newtheorem{lemma}[theorem]{Lemma}

\usepackage{babel}
\usepackage{authblk}
\usepackage{graphicx}
\usepackage{stmaryrd}
\usepackage{color}
\usepackage{adjustbox}
\usepackage{diagbox}

\usepackage{amsfonts} 
\definecolor{black}{rgb}{0,0,0}

\definecolor{red}{rgb}{1,0,0}

\definecolor{blue}{rgb}{0,0,1}

\makeatother

\newcommand{\cA}{\mathcal{A}}
\newcommand{\cE}{\mathcal{E}}

\newcommand{\cT}{\mathcal{T}}

\newcommand{\cN}{\mathcal{N}}

\newcommand{\bb}{\mathbf{b}}
\newcommand{\bc}{\mathbf{c}}
\newcommand{\bn}{\mathbf{n}}
\newcommand{\bx}{\mathbf{x}}
\newcommand{\bu}{\mathbf{u}}
\newcommand{\bv}{\mathbf{v}}

\newcommand{\bw}{\mathbf{w}}
\newcommand{\bbR}{\mathbb{R}}

\newcommand{\lam}{\lambda}
\newcommand{\la}{\langle}
\newcommand{\ra}{\rangle}

\newcommand{\pa}{\partial}

\graphicspath{{./figures/}}
\title{}

\title{\textbf{Physics-preserving IMPES based multiscale methods for immiscible two-phase flow in highly heterogeneous porous media} }
\author{Yiran Wang, Eric Chung and Shuyu Sun}

\begin{document}
	\maketitle
	\begin{abstract}
In this paper, we propose a physics-preserving multiscale method to solve an immiscible  two-phase flow problem, which is modeled as a coupling system consisting of Darcy's law and mass conservation equations. We use a new Physics-preserving IMplicit Pressure Explicit Saturation (P-IMPES) scheme in order to maintain the local conservation of mass for both phases. Besides, this scheme is unbiased and if the time step is smaller than a certain value, the saturation of both phases are bounds-preserving. When updating velocity, MGMsFEM serves as an efficient solver by computing the unknowns on a coarse grid. We follow the operation splitting techinque to deal with the two-phase flow. In particular, we use an upwind strategy to iterate the saturation explicitly and the MGMsFEM is utilized to compute velocity with a decoupled system on a coarse mesh. To show the efficiency and robustness of the proposed method, we design a set of interesting experiments. A rigorous analysis is also included to serve as a theoretical base of the method, which is well verified by the numerical results. Both simulations and analysis indicate that the method attains a good balance between accuracy and computation cost.
	\end{abstract}
\section{Introduction}

In scientific research and real applications, subsurface flow problems have drawn more and more attention \cite{mustapha2011efficient,fumagalli2011numerical,ginn1990inverse,forsyth1997monotonicity,wu1996consistent,panday2004fully,liang2007coupling} and a typical example is two-phase flow \cite{drew1983mathematical,ishii2010thermo,richter1983separated,ransom1984hyperbolic,spalding1981general} in porous media. Because of the high contrast and heterogeneity of the media, one needs to use a very fine grid to discretize the domain in order to obtain sufficient accuracy, which may result in much computation cost. Furthermore, in a time-dependent problem, one needs to solve the concerned equations at different times. If a single solve is computational extensive, the total cost is even larger. Hence, it is indispensable to develop some multiscale methods  to alleviate the computation burden. Upscaling \cite{chen2015new,gorguluarslan2015multilevel,boe1994analysis,chen2006adaptive} is a typical choice by averaging  quantities of interest on a coarse mesh. It can accelerate the computation but at the same time the accuracy suffers because of a loss of detailed information. Multiscale methods \cite{hou1997multiscale,chen2003mixed,hughes1998variational,arbogast2007multiscale,wheeler2012multiscale,cortinovis2014iterative,lunati2004multi}, on the other hand, can achieve a better balance between workload and precision. In particular, one needs to perform some fine-scale computation to generate some multiscale bases in the offline stage. Then, during each online process, simulations are carried out on a coarse mesh with the help of the bases. Since the bases contain some detailed information in the media, multiscale methods are more accurate than upscaling methods.

Among many multiscale methods, the Multiscale Finite Element Method (MsFEM) \cite{hou1997multiscale} is popular since it is straight-forward and relatively efficient. Solving a set of local problems on the fine grid, one can obtain a multiscale basis corresponding to each local region.  Based on MsFEM, a mixed formulation is further proposed to solve velocity and pressure in a decoupled system, which can satisfy the local conservation of mass and the method in this form is called the Mixed Multiscale Finite Element Method (MMsFEM) \cite{chen2003mixed,chung2015mixed}. Different from MsFEM, this method needs to construct velocity multiscale bases as well as pressure multiscale bases. Each velocity basis is supported in a local domain which is composed of two coarse elements sharing a common coarse edge. For each pressure basis, the support is a corresponding coarse element. However, since only one velocity and pressure bases can be obtained in each local region, the approximation effect  is not sufficiently good  especially when the media is highly heterogeneous. To this end, the Generalized Multiscale Finite Element Method (GMsFEM) \cite{eglp13,efendiev2013generalized,CHUNG2018419} is introduced to generate multiple bases in each region. The basis construction is divided into two steps. First, one needs to construct a local snapshot space which is spanned by a set of snapshot bases. Following this, some well-designed  spectral problems are utilized to reduce the dimension of snapshot space.  In particular, only  the eigenfunctions corresponding to smallest eigenvalues are included in the reduced space. In other words, only some representative modes are incorporated in the approximation space. However, it has been pointed out in \cite{efendiev2013generalized} that once the dimension of approximation space exceeds a particular value, the effect of including more bases is negligible. To further improve the accuracy, in \cite{chung2015residual,yang2018residual,efendiev2016online,wang2021online,wang2021comparison,wang2022local}, residual-driven bases are constructed. It has been shown in \cite{wang2021comparison} that the effect of using residual-driven bases is better than using more bases obtained from the eigenvalue problems.

Apart from computation efficiency, another concern is the local conservation of mass. A two-phase flow model is a coupling system composed of Darcy's law, conservation equations, constraints of saturations corresponding to two phases and the capillary pressure depending on the saturation corresponding to the wetting phase. To tackle the coupling system, some researches have been carried out including the standard IMplicit Pressure Explicit Saturation (IMPES) \cite{sheldon1959one,stone1961analysis} and some variations like HF-IMPES which is proposed by Hoteit and Firoozabadi \cite{hoteit2008numerical,hoteit2008efficient} to tackle the discontinuity of saturation. Both of two methods belong to IMplicit-EXplicit (IMEX) schemes but they can only promise the local conservation of mass for one phase, which is not sufficient in many situations. On the other hand, fully implicit method is unconditional stable but since it deals with all the unknown terms implicitly, one needs to endure a large computation cost. From the perspective of applications, it is not necessary to provide such a high accuracy. To this end, some new schemes \cite{chen2021new,chen2019fully1, chen2019fully2} are developed to satisfy the conservation for both phases. In \cite{chen2021new}, the authors introduce a new Physics-preserving IMPES (P-IMPES) scheme by calculating the total velocity as well as the velocity induced by the capillary pressure. During this process, an upwinding strategy is used for spatial discretization of the saturation. The scheme is proved in \cite{chen2021new} that it is local conservative for both phases, unbiased, and conditional bounds-preserving. In particular, the saturation of both phases are within their own bound when the time step is chosen less than a particular value.  

In this work, we propose a multiscale scheme based on the P-IMPES scheme. In particular, we utilize the multiscale mesh to update saturation and velocity separately. For saturation, we compute it on  a fine grid which is same as the one used in \cite{chen2021new}. On the other hand, to update velocity, since the heterogeneous  permeability field is incorporated which may result in much computation cost, we apply the MGMsFEM to improve the computation efficiency. We design a set of interesting experiments to show the error convergence as time step, mesh size and number of multiscale bases. We use two representative media, i.e., a highly heterogeneous media as well as a high-contrast media. The local conservation of mass is verified. Besides, errors of saturation decrease when a smaller time step or coarse mesh size is used. Moreover, we demonstrate the effect of enriching multiscale space. It is worth mentioning that when we enrich the approximation space, even a smaller number of residual-driven bases are more effective than the bases constructed from the spectral problems. Since we can adaptively choose different number of bases, this multiscale method is flexible. In the last section, we present a rigorous proof of the stability and convergence of the proposed method. The stability analysis includes that the local conservation property for both phases, unbiased property as well as the conditional bounds-preserving. Motivated from \cite{chen2021new},  we show that once the time partition is finer than one degree, then the stability can be obtained. Last but not least,  we derive the formulation of saturation error by finding the relation between saturation error as well as velocity error, which  serves as a firm base for the numerical performance of the scheme. It turns out that the simulation results are consistent to the final analysis conclusion.

The paper is composed of four parts. In the first section, we give the mathematical form of the concerned problem and present the P-IMPES scheme. Then, we state the MGMsFEM  and focus on the main topic, i.e., the reduced-order method based on IMPES. Numerical results are placed in the fourth section. Lastly, we prove the stability and convergence of this method.

	\section{Mathematical model and P-IMPES scheme}
\begin{figure}[!htbp] \centering
	\centering
	\includegraphics[width=0.9\textwidth]{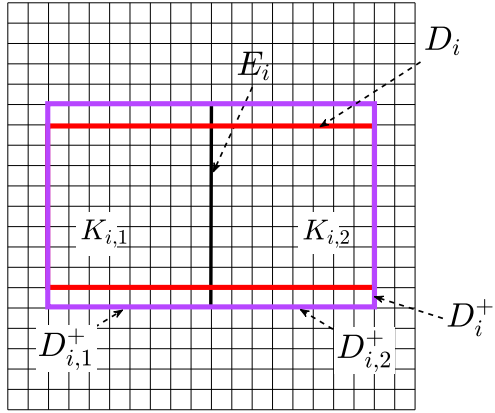}
	\caption{An illustration of a coarse-scale edge $E_i$, a corresponding coarse neighborhood $D_i$, and an oversampling neighborhood $D_i^+$.}
	\label{fig:mesh}
\end{figure}

	We first present the basic mathematical model for incompressible and immiscible two-phase flow in porous media. We denote the wetting and non-wetting phases by the subscripts $w$ and $n$, respectively. Our mathematical model is determined by utilizing the  conservation law, Darcy's law, the saturation  constraints and the capillary pressure. We consider a model with the gravity in porous media $\Omega\in \bbR^d(d=2,3)$ given as follows,
	\begin{eqnarray}
		\begin{aligned}
					\zeta\dfrac{\partial S_{\alpha}}{\partial t}+\nabla\cdot \bu_{\alpha}=q_{\alpha},  \quad&\text{ in } \Omega, \quad \alpha=w,n,\\
			\bu_{\alpha}=-\frac{k_{r\alpha}}{\mu_{\alpha}}\mathbf{K}(\nabla p_{\alpha}+\rho_{\alpha}g\nabla z), \quad&\text{ in } \Omega, \quad \alpha=w,n,\\
			S_n+S_w=1, \quad&\text{ in } \Omega,\\
			p_c(S_w)=p_n-p_w, \quad&\text{ in } \Omega.
			\label{model}
			\end{aligned}
		\end{eqnarray}
The initial and boundary conditions are specified below.
\begin{align*}
	S_{\alpha}&=S_{\alpha}^0,\quad t=0,\quad \alpha=w,n,\\
	p_{\alpha}&=p_{\alpha}^B,\quad \text{on }\Gamma_D,\quad \alpha=w,n,\\
	\bu_{\alpha}\cdot \bn&=g_{\alpha}^{N},\quad\text{on }\Gamma_N,\quad \alpha=w,n.
\end{align*}
	The mentioned notations are stated below.
	\begin{enumerate}
		\item Let $S_{\alpha}$, $\bu_{\alpha}$ and $p_{\alpha}$ be saturation, Darcy's velocity and pressure corresponding to phase $\alpha$. Total velocity is denoted by $\bu_t:=\bu_w+\bu_n$. 
	\item Define $\rho_\alpha$, $k_{r\alpha}$, $\mu_{\alpha}$ and  $q_{\alpha}$ as density, relative permeability, viscosity and sink/source term of phase $\alpha$.  Besides, $q_t=q_w+q_n$.   We define $\zeta$ as the porosity, which is a constant in $[0,1]$. Also, $\mathbf{K}$ is the absolute permeability tensor.
		\item Let $g$ be the magnitude of gravitational acceleration and $z$ be the depth.
		\item Phase mobility is denoted by  $\lambda_{\alpha}:=\frac{k_{r\alpha}}{\mu_{\alpha}}$. Total mobility is further defined by $\lambda_t=\lambda_w+\lambda_n$.
		\item We define fractional flow functions as $f_{w}=\lam_w/\lam_t$, $f_n=\lam_n/\lam_t$.
		\item We let the boundary be $\Gamma:=\partial \Omega$ and it can be partitioned as $\Gamma=\Gamma_D\cup \Gamma_N$, where $\Gamma_D$ and $\Gamma_N$ are Dirichlet and Neumann boundaries to solve $u_{\alpha}$. Moreover, to solve $S_{\alpha}$, we let $\Gamma=\Gamma_{\text{in}}\cup\Gamma_{\text{out}}$, where 
		$\Gamma_{\text{in}}=\{x\in \Gamma:\bu_t(\bx)\cdot \bn(\bx)<0\}$ is the inflow boundary,  $\Gamma_{\text{out}}=\{x\in \Gamma:\bu_t(\bx)\cdot \bn(\bx)\geq0\}$ is the outflow boundary, and $\bn$ is the unit outer normal vector to $\Gamma$.
		\item For phase $\alpha$, $p_{\alpha}^B$ is defined to be boundary condition of pressure on  $\Gamma_D$; $g_{\alpha}^N$ is set to be the  boundary condition of velocity  on  $\Gamma_N$; $S_{\alpha}^0$ is the initial saturation.
		\item Let $\xi_{\alpha}=\lambda_t\bw_{\alpha}$ with $\bw_{\alpha}=-\mathbf{K}(\nabla p_{\alpha}+\rho_{\alpha}g\nabla z)$, $\alpha=n,w$. Besides, we define $\xi_c=\xi_n-\xi_w.$

	\end{enumerate}
	Then we partition the global domain $\Omega$ into a coarse grid as well as a fine grid. Let $\cT_H$ be the coarse-scale  mesh.  Define $N_{e,c}$ and $N_{E,c}$ as the number of elements and edges in $\cT_H$. In particular, we let $\{K_i\}_{i=1}^{ N_{e,c}}$  be  coarse elements of   mesh size $H$, which form a conforming partition of $\Omega$, i.e. $\bar{\Omega}=\cup_{i=1}^{N_{e,c}}\{K_i\}$. We further let $\{E_i\}_{i=1}^{N_{E,c}}$ be the coarse edges. Besides, define $\cE_{c}(K_i)$ as the set of coarse edges of $K_i$ and $\cE_c:=\cup_{i=1}^{N_{e,c}}\cE_{c}(K_i)$. Among all the coarse-scale edges, we denote the number of all interior coarse-scale edges by $N_{\text{in},c}$. For each interior coarse edge, we can correspondingly define a local neighborhood $D_i$  composed of two coarse elements sharing the common coarse edge $E_i$. Since we use some local problems in the offline stages to capture fine-scale information in the media, we also need a fine grid. Correspondingly, we define $\cT_h$ as the fine-scale mesh with mesh size $h$.  For each $K_i$, we partition it into some fine-scale  elements with mesh size $h$. We define $\{K_i^h\}_{i=1}^{ N_{e,f}}$  as all fine-scale elements  with the total number $N_{e,f}$. Each $E_i$ is divided into some fine-scale edges with length $h$. We define $\cE_f(K_i)$ as all fine-scale edges on $\partial K_i$.  We further let $\cE_f=\cup_{i=1}^{N_{e,c}} \cE_{f}(K_i)$. We use $\{e_i\}_{i=1}^{N_{E,f}}$ to denote all the fine-scale edges, where $N_{E,f}$ is the total number. In the construction of multiscale bases, we  utilize the oversampling technique \cite{eglp13} to solve some local problems in oversampled local regions. In  Figure \ref{fig:mesh}, we use $D_i^+$ to denote an oversampled  local region corresponding to $D_i$. We remark that $D_i^+$ is generated by extending $D_i$ for  some fine-scale elements in one direction.

To compute the reference solutions, we solve \eqref{model} on a fine-scale mesh $\cT_h$. Here, we use the lowest Raviart-Thomas vector field (RT0) denoted by $V_h$ for fine-scale velocity bases. Besides, we define $Q_h$ as the pressure space spanned by piecewise constant functions and each basis corresponds to a particular fine-scale element.  And we use $V_h(D)$ and $Q_h(D)$ to represent  $V_h$ and $Q_h$ restricted in $D$. The P-IMPES scheme \cite{chen2021new} results the following equation system:
\begin{align}
	&\sum\limits_{\alpha} \beta_{\alpha}(\bu_t ^{h,n+1},q;S_w^{h,n})=(q_t,q).\label{scheme6}\\
	&(\zeta\dfrac{S_{\alpha}^{h,n+1}-S_{\alpha}^{h,n}}{t_{n+1}-t_n},q)+\beta_{\alpha}(\bu_t^{h,n+1},q;S_w^{h,n})=(q_{\alpha},q)+\sigma_{\alpha}\beta_c(\xi_c^{h,n+1},q;S_w^{h,n}), \quad \alpha=w \textbf{ or }n, \label{scheme1}\\
	&(\kappa_n^{-1}\bu_t^{h,n+1},\bv)-(p_w^{h,n+1},\nabla\cdot \bv)=(\kappa_n^{-1}f_n(S_w^{h,n})\xi_c^{h,n+1},\bv)-\int_{\Gamma_D} p_w^B \bv\cdot\bn-(\rho_wg\nabla z,\bv), \label{scheme2}\\
	&(\kappa_n^{-1}\xi_c^{h,n+1},\bv)=(p_c(S_w^{h,n}),\nabla\cdot \bv)-\int_{\Gamma_D} (p_n^B-p_w^B) \bv\cdot\bn-((\rho_n-\rho_w)g\nabla z,\bv).\label{scheme3}\\
	&(S_n^{h,n+1}+S_w^{h,n+1},q)=(1,q), \label{scheme4}\\
	&(p_n^{h,n+1}-p_w^{h,n+1},q)=(p_c(S_w^{h,n}),q). \label{scheme5}
\end{align}
The mentioned notations are defined below.
	\begin{enumerate}
		\item  Define $t_i$ as the $i$-th time step in a uniform partition of $[0,T]$, where $T$ is the final time. And $\Delta t=t_{n+1}-t_n$ for each $n$.
		\item Set $\sigma_w=1$ and $\sigma_n=-1$. 
		\item Let $\bu_t^{h,n}\in V_h(\Omega)$ and $\xi_c^{h,n}\in V_h(\Omega)$ be reference velocity solutions to \eqref{scheme3} and \eqref{scheme2} at time $t_n$. 
		\item  Use $S_{\alpha}^{h,n}\in Q_h$ to denote reference satuation of phase $\alpha$ at time $t_n$, which is the solution to \eqref{scheme1}.
		\item Define $p_{w}^{h,n}$ as reference pressure of wetting phase at time $t_n$ which is the solution to \eqref{scheme2}. Besides, $p_{n}^{h,n}$ is defined to be fine-grid pressure of non-wetting phase at time $t_n$ and it can be solved by \eqref{scheme5}.
		\item Define $\kappa_{n}:=\lambda_t(S_w^{h,n})\mathbf{K}$.
		\item If $q\in Q_h$ is piecewise constant, 
		\begin{align}
			\beta_{\alpha}(\bv,q;S_w^h)&=\sum\limits_{K^h\in \cT_h}\int_{\partial K^h} f_{\alpha}(S_{w,\alpha}^{*,h})\bv\cdot\bn q, \quad \alpha=w,n, \label{beta_alpha}\\
			\beta_{c}(\bv,q;S_w^h)&=\sum\limits_{K^h\in \cT_h}\int_{\partial K^h} f_{n}(S_{w,\alpha}^{*,h})f_{w}(S_{w,\alpha}^{*,h})\bv\cdot\bn q, \label{beta_c}	
		\end{align}
		where the upwind value $S_{w,\alpha}^{*,h}$  on $e\subset \pa K^h$ in the function $f_{\alpha}(S_{w,\alpha}^{*,h})$ is defined as follows:
		
		$S_{\alpha}^{*, h}|_e=\left\{\begin{array}{ll}\left.S_{\alpha}^{h}\right|_{K^h}, & \text { if }\left\{\boldsymbol{u}_{\alpha}^{h} \cdot \boldsymbol{n}\right\}_{e} \geq 0, \\ \left.S_{\alpha}^{h}\right|_{K^{h,1}}, & \text { if }\left\{\boldsymbol{u}_{\alpha}^{h} \cdot \boldsymbol{n}\right\}_{e}<0,\end{array} \quad S_{w, \alpha}^{*, h}= \begin{cases}S_{w}^{*, h}, & \alpha=w, \\ 1-S_{n}^{*, h}, & \alpha=n .\end{cases}\right.$
		
		Here $K^h\cap K^{h,1}=e$ and $\bn$ is a outward normal vector to $K^h$. For $e\subset \Gamma_{\text{in}}$, $S_{w,\alpha}^{*,h}|_e=S_w^{h}|_{K^h}$.
		We further define $\beta_t=\beta_w+\beta_n$.
	\end{enumerate}
	In the following, we use a matrix formulation.
	Before this, we assemble some matrices and vectors. At time $t_{k+1}$, we define
	\begin{align*}
		A^k[i,j]&=\int_{\Omega}\kappa_k^{-1} \bv_i\cdot\bv_j,\\
		A_n^k[i,j]&=\int_{\Omega}\kappa_k^{-1} f_n(S_w^{h,k})\bv_i\cdot\bv_j,\\
		C[i,j]&=\int_{\Omega}\nabla \cdot \bv_ip_j,\\
		P[i,j]&=\int_{\Omega}p_i p_j,\\
		D[j]&=\int_{\Gamma_D}\bv_j\cdot \bn,\\
		E[j]&=\int_{\Omega}g\nabla z\cdot \bv_j.\\
		F_{\alpha}[j]&=\int_{\Omega}q_{\alpha}p_j.\\
		P_c^k[j]&=p_c(S_w^{h,k})|_{K_j^h}.
	\end{align*}
Here, only $A^k$, $A_n^k$ and $P_c^k$ depend on time $t$ hence they should be updated as time. We further define  $F_t=F_w+F_n$. Besides, we define the corresponding matrices for $\beta_{\alpha}$ as 
\begin{align*}
	B_{\alpha}[i,j]=\beta_{\alpha}(\bv_j,p_i;S_w^{h,n}), \quad B_{t}[i,j]=\beta_{t}(\bv_j,p_i;S_w^{h,n}), \quad \alpha=w,n.
\end{align*}
Then, we have the following matrix formulation as follows:
	
	First, we solve $$A\xi_c^{h,n+1}=CP_c^n-(p_n^B-p_w^B)D-(\rho_n-\rho_w)E$$ to obtain $\xi_c^{h,n+1}$. This step corresponds to \eqref{scheme3}.
	
	Then we solve the following system 
	\begin{align*}
		\begin{bmatrix}
			A^n & -C\\
			B_t & 0
		\end{bmatrix}
		\begin{bmatrix}
			\bu_t^{h,n+1}\\ p_w^{h,n+1}
		\end{bmatrix}
		=\begin{bmatrix}
			A_n^n\xi_c-p_w^BD-\rho_wE\\ F_t
		\end{bmatrix},
	\end{align*}
which is associated with \eqref{scheme2} and \eqref{scheme6}.

	Finally, corresponding to \eqref{scheme1}, we can update saturation $S_w^{h,n+1}$ as follows,
	\begin{align*}
		 \zeta P\frac{S_w^{h,n+1}-S_w^{h,n}}{\Delta t} +B_w\bu_t^{h,n+1}=F_w+\sigma_{\alpha}B_c\xi_c^{h,n+1}.
	\end{align*}
	Hence,
	\begin{align*}
		S_w^{h,n+1}=S_w^{h,n}+\frac{1}{\zeta}\Delta t	 P^{-1}(F_w+\sigma_{\alpha}B_c\xi_c^{h,n+1}-B_w\bu_t^{h,n+1}).
	\end{align*}
\begin{table}[htbp!]
	\centering
	\begin{tabular}{c c}
		\hline\hline
		& P-IMPES Scheme\\
		\hline
		Step 1: &Seek $p_{c}^{h,n+1}=p_n^{h,n+1}-p_{w}^{h,n+1}\in Q_h$ and $\xi_c^{h,n+1}\in U_h$ by \eqref{scheme5} and \eqref{scheme3}.\\
		\hline
		Step 2: & Use \eqref{scheme6} and \eqref{scheme2} to solve $p_w^{h,n+1}$ and $\bu_t^{h,n+1}$. \\
		&Then $p_n^{h,n+1}$ can be updated by \\
		&$p_n^{h,n+1}=p_{nw}^{h,n+1}+p_w^{h,n+1}.$\\
		\hline
		Step 3: & Update the wetting phase saturation $S_w^{h,n+1}$ by \eqref{scheme1} for $\alpha=w$. \\
		&Then the non-wetting phase saturation is updated by \eqref{scheme4}.\\
		\hline \hline
	\end{tabular}
\end{table}
	\section{Mixed generalized multiscale finite element method}
	
	We use a mixed formulation of the generalized multiscale finite element method (MGMsFEM) \cite{chung2015mixed}. The MGMsFEM contains two stages of offline construction, denoted by the offline stage I and offline stage II. In the first stage, we first construct a snapshot space by solving a set of local problems. Using a set of well-designed spectral problems, we can further reduce the snapshot space into a smaller multiscale space. Based on the multiscale bases, we can solve the previous problem on a coarse grid, which is much more computationally efficient. The second stage is constructing some residual-driven bases. We perform the construction of these  bases in each oversampled local region. The final multiscale space is linearly spanned by the multiscale bases constructed in the above two stages.
	
	We remark that the construction of multiscale bases in two stages are designed by different motivations. For the first stage, since we have constructed a large approximation space, i.e. the snapshot space, we need to reduce it properly by extracting the most important modes from the snapshot space. To this end, we solve a set of spectral problems and only incorporate  a part of eigenfunctions in the multiscale space in the offline stage I. On the other hand, the goal of second stage is to enhance the approximation accuracy of multiscale space. We utilize the residuals of solving the coarse-scale system corresponding to previous multiscale space to construct new bases, which will be incorporated into the multiscale space. We  present the detailed process of construction as follows.
	
	Let $\kappa_0:=\lambda_t(S_w^0) \mathbf{K}$, which is used  to construct the multiscale bases. We use traditional notations in Sobolev space such as $L^2(\Omega), H^1(\Omega), H_0^1(\Omega)$. Besides, we use  $H(\text{div};\Omega;\kappa_0^{-1})$ to denote the Sobolev space containing vector fields $\bv$ with 
	$\bv\in L^2(\Omega)^2$ and $\nabla\cdot \bv\in L^2(\Omega)$. The corresponding norm of $H(\text{div};\Omega;\kappa^{-1})$ is $\|\bv\|_{H(\text{div};\Omega;\kappa^{-1})}$. We define the mentioned notations of inner products and norms as follows. For all $\bu,\bv\in (L^2(\Omega))^2$, $p,q\in H_0^1(\Omega)$, we let
	\begin{align*}
		\la p,q\ra_{\Omega}&=\int_{\Omega} pq, \quad \cA_{\Omega}(\bu,\bv)=\int_{\Omega} \kappa_0^{-1} \bu \cdot \bv, \\
		\|p\|_{L^2(\Omega)}^2&=\la p,p\ra_{\Omega}, \quad \|\bu\|_{\kappa_0^{-1},\Omega}^2=\cA_{\Omega}(\bu,\bu),\quad 
		\|\bu\|_{H(\text{div};\Omega;\kappa_0^{-1})}^2=\|\bu\|_{\kappa_0^{-1},\Omega}^2+\|\nabla \cdot \bu\|_{L^2(\Omega)}^2.
	\end{align*}

 To perform the computation on $\cT_H$, we define $Q_H$ to be the pressure space composed of piecewise contant functions in the coarse grid mesh $\cT_H$. In the following, we show the process of computing multiscale bases in two stages.
	\subsection{Offline stage I}
	First of all, we construct a snapshot space using the following local problems. For each $i$, we define $L_i$ to be the number of local snapshots in $D_i$. We compute $\psi_i^j$ and $p_i^j$ for $j=1,\ldots, L_i$ in $K_{i,1}$ and $K_{i,2}$ separately, where $K_{i,1}\cup K_{i,2}=D_i$ and $K_{i,1}\cap K_{i,2}=E_i$. In particular, we solve
	\begin{eqnarray}
		\begin{aligned}
			\kappa_0^{-1} \psi_i^j+\nabla p_i^j&=0 \text{ in } K_{i,p},\\
			\nabla\cdot\psi_i^j&=\alpha_i^p \text{ in }K_{i,p},\\
			\psi_i^j\cdot \bn_i&=0 \text{ on }\partial D_i\cap\partial K_{i,p},\\
			\psi_i^j\cdot \bn_i&=\delta_i^j\text{ on }E_i,\label{snap1}		
		\end{aligned}
	\end{eqnarray}
	for $p=1,2$. We set $\alpha_i^p$ to be a constant in $K_{i,p}$, which satisfies $\alpha_i^1+\alpha_i^2=0$. Moreover, $\delta_i^j$ is defined by 
	\begin{equation*}
		\delta_i^j=\left\{
		\begin{array}{cc}
			1,&\text{ on } e_j,\\
			0,& \text{ on }E_i\setminus e_j,
		\end{array}
		j=1,\ldots,L_i.
		\right.
	\end{equation*}
	We need to emphasize that $\alpha_i^j$ and $\delta_i^j$ should satisfy a compatible condition that $\int_{K_{i,p}} \alpha_i^j=\int_{E_i} \delta_i^j$ for $p=1,2$.
	Instead of solving \eqref{snap1} directly, we solve $\psi_i^j$ and $p_i^j$ with a variational formulation as follows. For all $\bu\in V_h(K_{i,p})$, $q\in Q_h(K_{i,p})$,
	\begin{align}
		\cA_{K_{i,p}} (\psi_i^j,\bu)-\la p_i^j,\nabla\cdot\bu\ra_{K_{i,p}}&=0,\\
		\la\nabla\cdot\psi_i^j,q\ra_{K_{i,p}} &=\la \alpha_i^j,q \ra_{K_{i,p}},
		\label{snap2}
	\end{align}
	where the boundary conditions and  the compatible requirements for $\alpha_i^j$ and $\delta_i^j$ are the same as the one used in \eqref{snap1}. We remark that each $\psi_i^j$ is further extended to the global domain $\Omega$ by defining $\psi_i^j=0$ outside of $D_i$. For simplicity, we use the same notation for the global version.
	We let $V_{\text{snap}}^i:=\text{snap}\{\psi_i^1,\ldots, \psi_i^{L_i}\}.$ Then the final snapsot space is defined  by $V_{\text{snap}}:=\bigoplus_{i=1}^{N_{\text{in},c}} V_{\text{snap}}^i$. However, since $V_{\text{snap}}$ contains $O(1/h)$ snapshots, we  use a set of local spectral problems to reduce it to a smaller approximation space.
	
	We first introduce the following two bilinear inner products. In each $D_i$ and $E_i$, we define
	\begin{eqnarray}
		\begin{aligned}
			a_i(\bv,\bu)&=\int_{E_i}\kappa_0^{-1}(\bv\cdot \bn_i)(\bu\cdot \bn_i),\\
			s_i(\bv,\bu)&=\frac{1}{H}\left(\int_{D_i}\kappa_0^{-1}\bv\cdot \bu+\int_{D_i}(\nabla\cdot\bv)(\nabla\cdot\bu)\right), \label{loc_spe}
		\end{aligned}
	\end{eqnarray}
	for all $\bv, \bu\in V_{\text{snap}}^i$.
	We then solve 
	\begin{align}
		a_i(\bv_i^j,\bu)=\lambda_i^j s_i(\bv_i^j,\bu), \quad \forall \bu\in V_{\text{snap}}^i.
	\end{align} 
	We arrange $\lambda_i^j$ in an ascending order and we include the eigenvectors corresponding to the smallest eigenvalues in the local multiscale space.  In particular, we define $V_{\text{ms},0}^i=\{\bv_i^1,\ldots, \bv_i^{l_i}\}$, where $l_i$ is defined as the number of bases in $V_{\text{ms},0}^i$.  We finally define $V_{\text{ms},0}=\bigoplus_{i=1}^{N_{\text{in},c}}V_{\text{ms},0}^i$. Here, the subindex  $0$ means no bases enrichment is performed.
	
	We remark that the design of local spectral problems \eqref{loc_spe} to $V_{\text{ms},0}$ is motivated from analysis. However, one can observe through out the whole process in the first stage, no information of source term $f_{\alpha}$, $\alpha=w,n$, is taken into account, which will restrict the approximation accuracy to some extent. To this end, we also need a second offline stage to construct some new bases containing some useful information not included in  $V_{\text{ms},0}$.
	\subsection{Offline stage II}
	In this part, we will perform bases enrichment to the multiscale space obtained in Offline stage I. The construction at this stage is performed by iterations and we define $V_{\text{ms},k}$ to be the multiscale space after $k$ iterations. In this way, $V_{\text{ms},0}$ is the result of Offline stage I. We use $(\bv_{\text{ms}}^k,p_{\text{ms}}^k)$ to denote the solutions of \eqref{ms_k} sought in $V_{\text{ms},k}$ and $Q_H$. Then we define the residual operator. Suppose $D\subset \Omega$ is a domain. We define $V_D:=\cup_{D_i\subset D} V_{\text{snap}}^i$ and $\hat{V}_D$ is the corresponding divergence free subspace of $V_D$.  We then state the procedures of stage II as follows. At iteration $k+1$,
	the first step is to solve $(\bv_{\text{ms}}^{k},p_{\text{ms}}^{k})$ by
	\begin{eqnarray}
		\begin{aligned}
			\int_{\Omega}\kappa_0^{-1}\bv_{\text{ms}}^k\cdot \bu-\int_{\Omega}\nabla\cdot\bu p_{\text{ms}}^k&=0,\quad \forall \bu\in V_{\text{ms,k}},\\
			\int_{\Omega}\nabla\cdot\bv_{\text{ms}}^k q&=\int_{\Omega}q_tq,\quad \forall q\in Q_{H},
			\label{ms_k}
		\end{aligned}
	\end{eqnarray}
	where $q_t=q_w+q_n$.
	Before stepping into the second step, we introduce a residual operator. Given a pair of multiscale solution $(\bv_{\text{ms},k},p_{\text{ms},k})$, the residual operator $R_D^{(k)}$ on $V_D$ is defined by 
	\begin{align}
		R_D^{(k)}(\bu)=\int_{D} \kappa_0^{-1}\bv_{\text{ms},k}\bu-\int_{D} \nabla\cdot\bu p_{\text{ms},k}, \quad \forall \bu\in V_{D}.\label{loc_res}
	\end{align}
	Restricted on $\hat{V}_D$, $R_{D}^{(k)}$ can be simplified as 
	\begin{align*}
		R_{D}^{(k)}(\bu)=\int_{D}\kappa_0^{-1} \bv_{\text{ms},k}\cdot \bu, \quad \forall \bu\in \tilde{V}_{D}.
	\end{align*}
	We denote the operator norm of the residual $R_D^{(k)}$ by $\|R_D^{(k)}\|$.
	The second step is to compute the global residual $\|R_{\Omega}^{(k)}\|$ and compare it with a given tolerance $\tau$. If $\|R_{\Omega}^{(k)}\|>\tau$, we perform the following three steps. Otherwise, we can terminate the enrichment.
	
	The third step is to select a set of local regions to perform bases enrichment. Suppose $D_1,\ldots, D_{B_k}$ form a nonoverlapping partition of $\Omega$, where $B_k$ is defined to be the number of local regions partitioning $\Omega$. For each $D_i$, we generate an oversampled region $D_i^+$ by enlarging $D_i$ for some fine-scale elements. In each $D_i^+$, we compute $R_{D_i^+}^{(k)}$ as defined in \eqref{loc_res}.
	
	Then, we perform the construction of new bases as below. In each $D_i^+$, we compute a corresponding $\eta_i$ by solving
	\begin{align*}
		\cA_{D_i^+}(\eta_i, \bu)=R_{D_i^+}^{(k)}(\bu), \quad \forall \bu \in V_{D_i^+}.
	\end{align*}
	We remark that using this construction, each $\eta_i$ is supported in $D_i^+$. Each $D_i$ and $D_i^+$ corresponds to an inner coarse edge $E_i$, which is showed in the Figure \ref{fig:mesh}. We restrict $\eta_i\cdot \bn_i$ on $E_i$ and normalize it to get $\lambda_i$.
	
	The final step is that we use the boundary information of the solution in the last step and solve the following problem in each $K_{i,p}$ individually with $i=1,\ldots,B_K$ and $p=1,2$,
	\begin{eqnarray}
		\begin{aligned}
			\kappa_0^{-1}\phi_i+\nabla p_i&=0\text{ in } K_{i,p},\\
			\nabla\cdot\phi_i&=\gamma_i^p\text{ in } K_{i,p},\\
			\phi_i\cdot \bn_i&=\lambda_i\text{ on }E_i,\\
			\phi_i\cdot\bn_i&=0\text{ on } \partial D_i\cap\partial K_{i,p},\label{online}
		\end{aligned}
	\end{eqnarray}
	where $\gamma_i^p$ is a constant in $K_{i,p}$ such that $\gamma_i^1+\gamma_i^2=0$. Besides, $\gamma_i^p$ is chosen to satisfy the compatible condition that $\int_{K_{i,p}} \gamma_i=\int_{E_i} \lambda_i,$ for both $K_{i,1}$ and $K_{i,2}$ satisfying  $K_{i,1}\cup K_{I,2}=D_i$. Similar to the construction of $V_{\text{snap}}$, we extend each $\phi_i$ to $\Omega$ by setting $\phi_i=0$ elsewhere. We keep using $\phi_i$ for simplicity of notations. Finally, we include the new bases in the previous multiscale space, i.e. $V_{\text{ms},k+1}^{i}=V_{\text{ms},k}^{i}\oplus \text{span}\{\phi_i\}$. Then $V_{\text{ms},k+1}=\bigoplus_{i=1}^{N_{\text{in},c}}V_{\text{ms},k+1}^{i}$.
	\subsection{Combine MGMsFEM with P-IMPES}
	From previous MGMsFEM, one has obtained a multiscale space for velocity, $V_{\text{ms}}=\{\phi_1,\ldots, \phi_{N_{\text{ms}}}\}$, where $N_{\text{ms}}$ is the number of multiscale bases in the final multiscale space. We recall that $Q_H=\{q_1^H,\ldots, q_{N_{e,c}}^H\}$ is a piecewise constant space corresponding to the coarse-scale mesh $\cT_H$. In particular, each $q_i^H$ is constant $1$ on the corresponding coarse element and vanishes outsides. We recall that $\xi_c$ and $\bu_t$ are solutions to \eqref{scheme3} and \eqref{scheme2}, respectively. Thus, one can use $V_{\text{ms}}$ and $Q_H$ to compute the pressure equation system in a coarse scale. As for the transport part, one can use the mixed finite element method to update satuation, which is the same with the P-IMPES scheme. Hence, the proposed alogrithm is a combination of MGMsFEM and P-IMPES, which can be denoted as MS-P-IMPES. We solve the following system: $\forall \bv_H\in V_{\text{ms}}$, $q_H\in Q_H$, and $q\in Q_h$,
	\begin{align}
		&\sum\limits_{\alpha}\beta_{\alpha}(\bu_t ^{H,n+1},q_H;S_w^{H,n})=(\tilde{q}_t,q_H).\label{ms-scheme1}\\
		&(\zeta\dfrac{S_{\alpha}^{H,n+1}-S_{\alpha}^{H,n}}{t_{n+1}-t_n},q)+\beta_{\alpha}(\bu_t^{H,n+1},q;S_w^{H,n})=(q_{\alpha},q)+\sigma_{\alpha}\beta_c(\xi_{c}^{H,n+1},q;S_w^{H,n}), \quad \alpha=w \textbf{ or }n, \label{ms-scheme2}\\
		&((\kappa_n^H)^{-1}\bu_t^{H,n+1},\bv_H)-(p_w^{H,n+1},\nabla\cdot \bv_H)=((\kappa_n^H)^{-1}f_n(S_w^{H,n})\xi_c^{H,n+1},\bv_H)-\int_{\Gamma_D} p_w^B \bv_H\cdot\bn-(\rho_wg\nabla z,\bv_H), \label{ms-scheme3}\\
		&((\kappa_n^H)^{-1}\xi_c^{H,n+1},\bv_H)=(\tilde{p}_c(S_w^{H,n}),\nabla\cdot \bv_H)-\int_{\Gamma_D} (p_n^B-p_w^B) \bv_H\cdot\bn-((\rho_n-\rho_w)g\nabla z,\bv_H).\label{ms-scheme4}\\
		&(S_n^{H,n+1}+S_w^{H,n+1},q)=(1,q), \label{ms-scheme5}\\
		&(p_n^{H,n+1}-p_w^{H,n+1},q_H)=(\tilde{p}_c(S_w^{H,n}),q_H). \label{ms-scheme6}
	\end{align}
We define the mentioned notations below.
\begin{enumerate}
	\item Let $\bu_t^{H,n}\in V_{\text{ms}}$, $\xi_c^{H,n}\in V_{\text{ms}}$ be multiscale solutions to \eqref{ms-scheme3} and \eqref{ms-scheme4}\ at time $t_n$.
	\item  Define $p_{w}^{H,n}\in Q_H$ as coarse-grid pressure of wetting phase at time $t_n$, which is the solution to \eqref{ms-scheme3}. Besides, we define $p_{n}^{h,n}$ as fine-grid pressure of non-wetting phase at time $t_n$ and it can be solved by \eqref{ms-scheme6}.
	\item Let $S_{\alpha}^{H,n}\in Q_h$ be multiscale saturation of phase $\alpha$ at time $t_n$, which is the solution to \eqref{ms-scheme2} based on $\bu_t^{H,n}\in V_{\text{ms}}$ and $\xi_c^{H,n}\in V_{\text{ms}}$.  
	\item Let $\tilde{p}_c(S_{\alpha}^{H,n})\in Q_H$ be an averaged capillary pressure of ${p}_c(S_{\alpha}^{H,n})$ in each coarse element. In particular, $\tilde{p}_c(S_{\alpha}^{H,n})|_K=\int_{K}\tilde{p}_c(S_{\alpha}^{H,n})/|K|$, $\forall K\in \cT_H$.
	\item Let $\kappa_{n}:=\lambda_t(S_w^{h,n})\kappa_{0}$ and $\kappa_{n}^H:=\lambda_t(S_w^{H,n})\kappa_{0}$.
	\item Let $\kappa_{\text{min},n}:=\min\limits_{\bx\in \Omega} \kappa_{n}(\bx)$ and $\kappa_{\text{max},n}:=\max\limits_{\bx\in \Omega} \kappa_n(\bx)$.
	Correspondingly, 
	
	$\kappa_{\text{min},n}^H:=\min\limits_{\bx\in \Omega} \kappa_{n}^H(\bx)$ and $\kappa_{\text{max},n}^H:=\max\limits_{\bx\in \Omega} \kappa_n^H(\bx)$.
\end{enumerate}

For a better illustration, we also use a matrix presentaion. Before this, we then define matrix formulation for the multiscale space. 
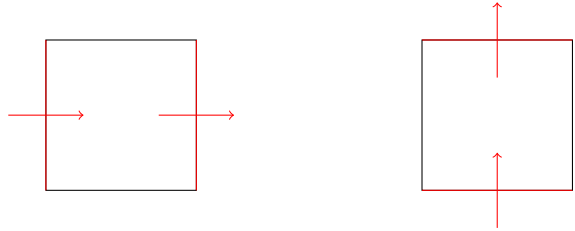
\begin{figure}[!htbp] \centering
	\centering

\begin{tikzpicture}

	\draw[black] (0,0) rectangle (2,2);
    \draw[red] (2,0)--(2,2);
    \draw[red] (0,0)--(0,2);
	\draw[->][red] (1.5,1)--(2.5,1);
	\draw[->][red] (-0.5,1)--(0.5,1);
	
	\draw[black] (5,0) rectangle (7,2);
	\draw[red] (5,2)--(7,2);
	\draw[red] (5,0)--(7,0);
	\draw[->][red] (6,1.5)--(6,2.5);
	\draw[->][red] (6,-0.5)--(6,0.5);
\end{tikzpicture} 
\caption{ An illustration of the relationship of edges and correponding normal vectors.}
\label{normal vector}
\end{figure}

For each velocity multiscale basis $\phi_i$, we evaluate it on the fine-scale edges. In particular, let $\bb_i\in \bbR^{N_{E,f}}$ be a vector corresponding to $\phi_i$  defined by  $\bb_i[j]=\phi_i\cdot \bn_j$, where $\bn_j$ is a unit normal vector to the fine-scale edge $e_i$ with pointing to the right or upward as displayed in Figure \ref{normal vector}.
Then we define $\Phi_v\in \bbR^{N_{E,f}\times N_{\text{ms}}}$ as $\Phi_v=[\bb_1,\ldots, \bb_{N_{\text{ms}}}]$. We also need to define a matrix for coarse-scale pressure space $Q_H$. For each $q_i^H$, we define a corresponding vector $\bc_i\in \bbR^{N_{e,f}}$ such that $\bc_i[j]=1$ for $K_j^h\subset K_i$ and $\bc[j]=0$ elsewhere.  Let $\Phi_p\in \bbR^{N_{e,f}\times N_{e,c} }$ be the matrix $Q_H$ and $\Phi_p:=[\bc_1,\ldots,\bc_{N_{e,c}}]$.
Furthermore, we can generate the following matrices by relating them with the previous matrices asembled on the fine-scale mesh. At time $t_{k+1}$,
	\begin{align*}
	\tilde{A}^{k}[i,j]&=\int_{\Omega}(\kappa_k^H)^{-1}\phi_i\cdot \phi_j=\Phi_v^TA^k\Phi_v,\\
	\tilde{A}_n^k[i,j]&=\int_{\Omega}(\kappa_k^H)^{-1} f_n(S_w^{H,k})\phi_i\cdot \phi_j=\Phi_v^TA_n^k\Phi_v,\\
	\tilde{C}[i,j]&=\int_{\Omega}\nabla \cdot \phi_i\cdot q_j^H=\Phi_v^TC\Phi_p,\\
	\tilde{P}[i,j]&=\int_{\Omega}q_i^Hq_j^H=\Phi_p^TP\Phi_p,\\
	\tilde{D}[j]&=\int_{\Gamma_D}\phi_j\cdot \bn=\Phi_v^TD,\\
	\tilde{E}[j]&=\int_{\Omega}g\nabla z\cdot \phi_j=\Phi_v^TE.\\
	\tilde{F}_{\alpha}[j]&=\int_{\Omega}q_{\alpha}q_j^H=\Phi_p^TF.\\
	\tilde{P}_c^k[j]&=\tilde{p}_c(S_w^{H,k})|_{K_j}.
\end{align*}
Correspondingly, only $\tilde{A}^k$, $\tilde{A}_n^k$ and $\tilde{P}_c^k$ should be updated as time. We define  $\tilde{F}_t=\tilde{F}_w+\tilde{F}_n$. Besides, we define $\tilde{B}_{\alpha}$ as $\tilde{B}_{\alpha}=\Phi_p^TB_{\alpha}\Phi_v$.  Also, $\tilde{B}_t=\tilde{B}_w+\tilde{B}_n$. Based on the above matrices, we can proceed to the coarse-scale computations. Thus, we obtain reduced-order solutions, which are indeed coefficients with respect to bases. Here, we define $\tilde{\xi}_c^{H,n+1}$, $\tilde{\bu}_t^{H,n+1}$ and $\tilde{p}_w^{H,n+1}$ to be reduced solutions corresponding to  $\xi_c^{H,n+1}$, $\bu_t^{H,n+1}$ and $p_w^{H,n+1}$. In particular, their relationships are characterized as follows:
\begin{eqnarray}
	\begin{aligned}
		\xi_c^{H,n+1}&=\sum_{i=1}^{N_{\text{ms}}} \tilde{\xi}_c^{H,n+1}[i]\phi_i,\\
		\bu_t^{H,n+1}&=\sum_{i=1}^{N_{\text{ms}}} \tilde{\bu}_t^{H,n+1}[i]\phi_i,\\
		p_w^{H,n+1}&=\sum_{i=1}^{N_{e,c}} \tilde{p}_w^{H,n+1}[i]q_i^H.
		\label{coarse_fine}
	\end{aligned}
\end{eqnarray}

Corresponding to \eqref{ms-scheme4}, the first step is  solving  $$\tilde{A}^n\tilde{\xi}_c^{H,n+1}=\Phi_v^TC\tilde{P}_c-(p_n^B-p_w^B)\tilde{D}-(\rho_n-\rho_w)\tilde{E}$$ to obtain $\tilde{\xi}_c^{H,n+1}$.

Then,  from \eqref{ms-scheme1} and \eqref{ms-scheme3}, we compute the following system 
\begin{align*}
	\begin{bmatrix}
		\tilde{A}^n & -\tilde{C}\\
		\tilde{B}_t & 0
	\end{bmatrix}
	\begin{bmatrix}
		\tilde{\bu}_t^{H,n+1}\\\tilde{p}_w^{H,n+1}
	\end{bmatrix}
	=\begin{bmatrix}
		\tilde{A}_n^n\tilde{\xi}_c^{H,n+1}-p_w^B\tilde{D}-\rho_w\tilde{E}\\ \tilde{F}_t
	\end{bmatrix}.
\end{align*}
Then $\bu_t^{H,n+1}$ and $p_w^{H,n+1}$ are obtained using  \eqref{coarse_fine}.

To satisfy the mass conservation property, we use a postprocessing techique. We define $\bu_{t,p}^{H,n+1}$ to be the velocity solution solved by postpocessing. Let $S_{q_t}$ be the set of coarse elements where $q_t$ is not a constant. For each $K\in S_{q_t}$, we solve a local problem on the fine grid as follows. 

For all $\bv\in V_{h}(K)$ and $q\in Q_h(K)$,
\begin{eqnarray*}
	\begin{aligned}
		\int_{K}\kappa^{-1}\bu_{t,p}^{H,n+1}\cdot\bv+\int_{K} \nabla p_w^{H,n+1}\cdot\bv&=0, \\
\int_{K}\nabla\cdot \bu_{t,p}^{H,n+1} q&=\int_{K} q_t q,\\
\bu_{t,p}^{H,n+1}&=\bu_{t}^{H,n+1},\text{ on }\pa K.	\label{postpro}	
	\end{aligned}
\end{eqnarray*}

We adjust the final velocity solution by $\bu_{p}^{H,n+1}:=\bu_{t,p}^{H,n+1}$, $\forall K\in S$. Based on the postpocessing procedure, \eqref{ms-scheme1} is equivalent to 
\begin{align}
	\sum\limits_{\alpha}\beta_{\alpha}(\bu_t ^{H,n+1},q;S_w^{H,n})=(q_t,q),\quad \forall q\in Q_h. \label{pp}
\end{align}
Indeed, from the definition of $\alpha_i^p$ in \eqref{snap1}, we can know $\nabla\cdot \psi_i^j$ is a constant in each coarse element for all $i$ and $j$, which results that $\nabla \cdot \bu_t ^{H,n+1}$ is a piesewise constant function. For $K\notin S$, $q_t$ is a constant hence \eqref{pp} and\eqref{ms-scheme1} are equivalent. At the end, we conclude the MS-P-IMPES Scheme in the Table \ref{table_ms-p-impes} below.
\begin{table}[htbp!]
	\centering
	\begin{tabular}{c c}
		\hline\hline
		& MS-P-IMPES Scheme\\
		\hline
		Step 1: &Seek $p_{n,w}^{H,n+1}=p_n^{H,n+1}-p_{w}^{H,n+1}\in Q_H$ and $\tilde{\xi}_c^{H,n+1}\in V_{\text{ms}}$ by \eqref{ms-scheme6} and \eqref{ms-scheme4}.\\
		\hline
		Step 2: & Use \eqref{ms-scheme1} and \eqref{ms-scheme3} with postprocessing to solve $p_w^{h,n+1}$ and $\bu_t^{H,n+1}$. \\
		&Then $p_n^{H,n+1}$ can be updated by \\
		&$p_n^{H,n+1}=p_{nw}^{H,n+1}+p_w^{H,n+1}.$\\
		\hline
		Step 3: & Update the wetting phase saturation $S_w^{H,n+1}$ by \eqref{ms-scheme2} for $\alpha=w$. \\
		&Then the non-wetting phase saturation is updated by \eqref{ms-scheme5}.\\
		\hline		\hline
	\end{tabular}
\caption{MS-P-IMPES Scheme: a combination of MMsFEM and P-IMPES. We use postprocessing in the coarse elements with nonconstant source.}
\label{table_ms-p-impes}
\end{table}

	\section{Numerical simulations}
	In this section, we show some numerical examples to verify the performance of MS-P-IMPES scheme. First of all, we verify that the conservation of mass is satisfied by the proposed scheme. Furthermore, we show the convergence of error as the decreasing of time step size, mesh size and increasing the number of multiscale bases. Besides, we demonstrate the approximation effects in capturing fine-scale information through presenting solutions in some local regions. Last but not least, we explore the effect of source terms on saturation dynamics. 
	
\begin{figure}[!htbp] \centering
		\centering
		\subfigure[$\ln\kappa_1(x)$]
		{ \includegraphics[width=0.45\textwidth]{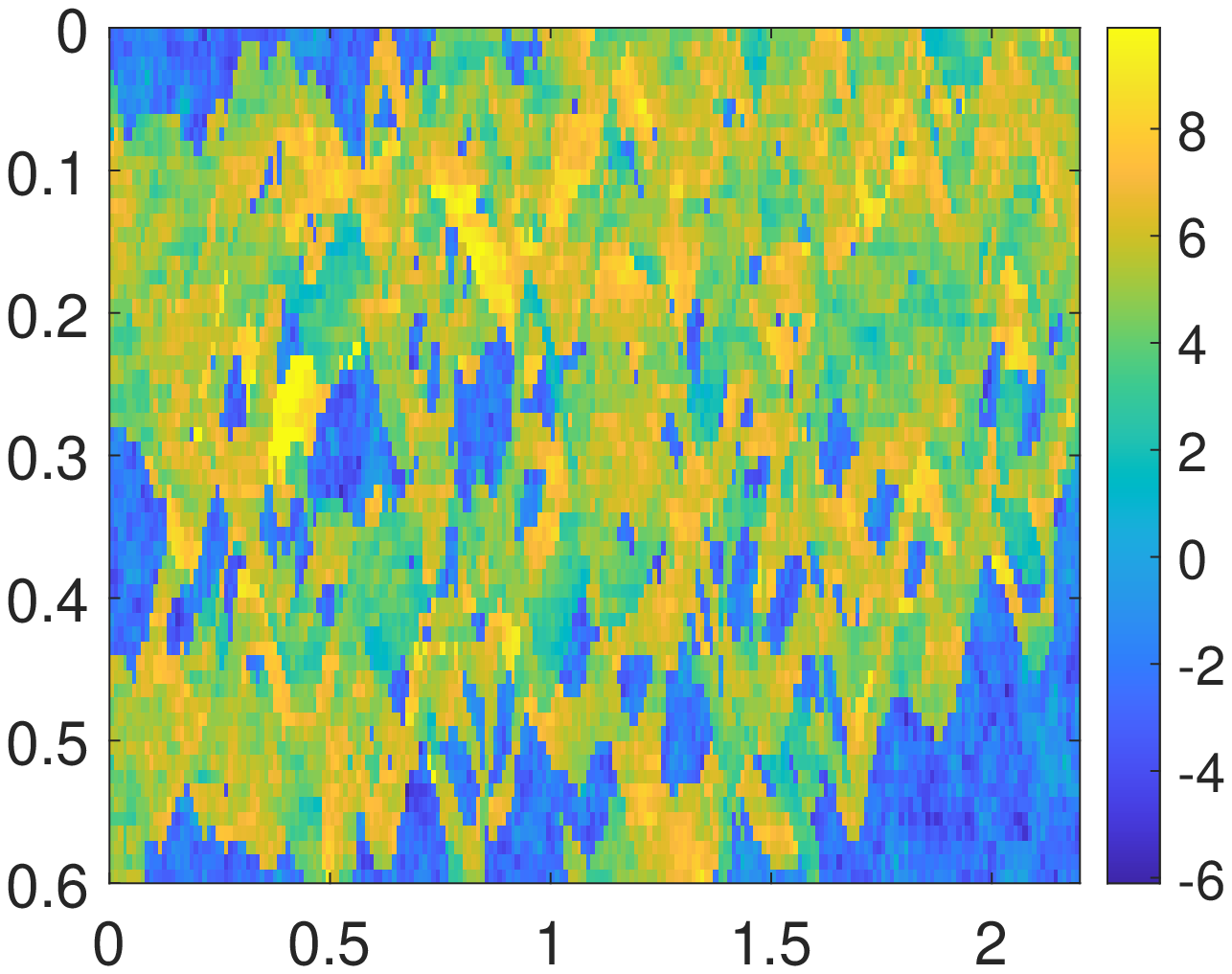}}
		\subfigure[$\ln\kappa_2(x)$]
		{ \includegraphics[width=0.45\textwidth]{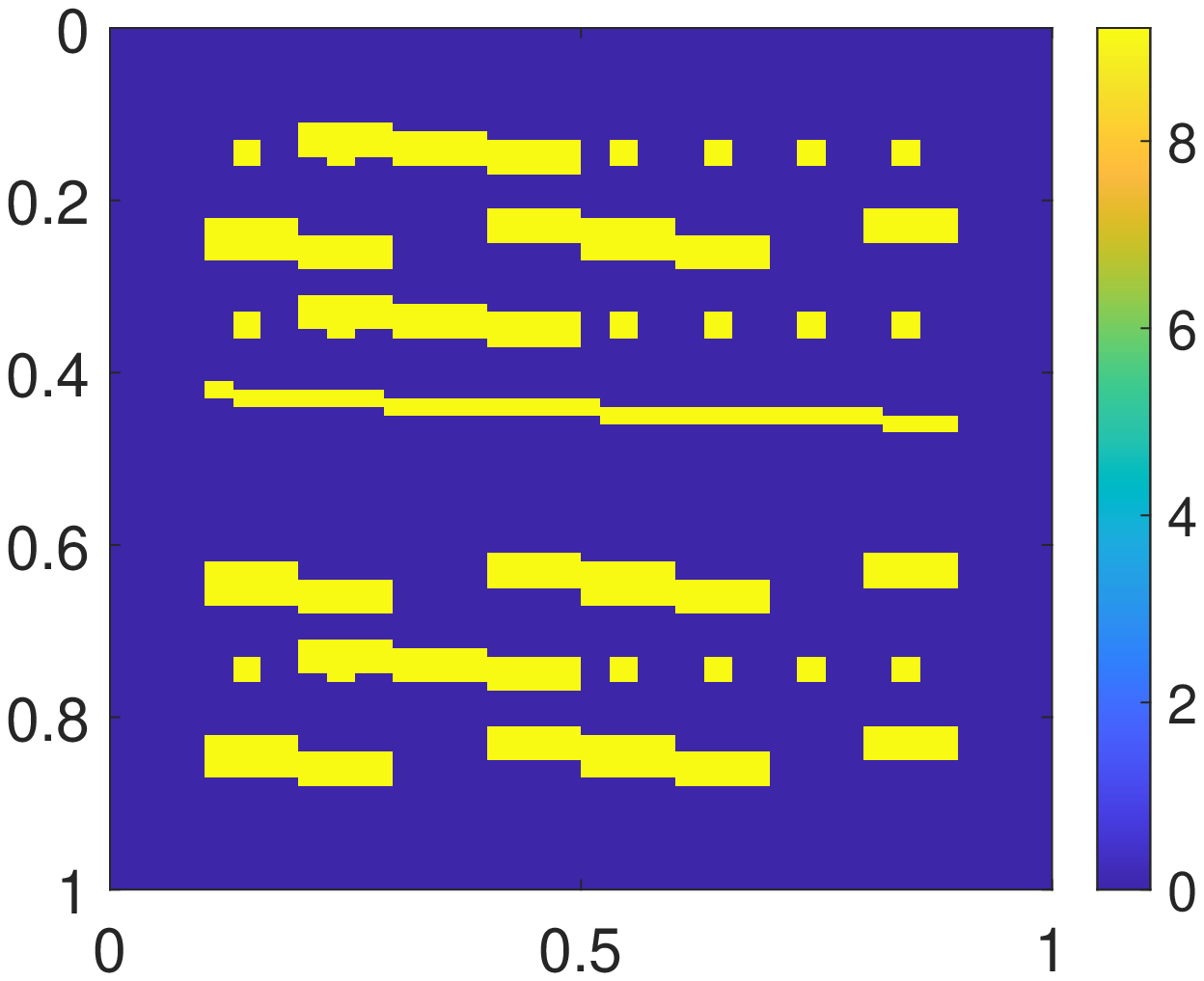}}
		\caption{Examples of heterogeneous permeability fields in a natural log scale}\label{medium}
	\end{figure}

	Through the numerical simulations, we use the SPE10 benchmark model as well as a high-contrast model, denoted by $\kappa_1$ and $\kappa_2$, whose natural logarithms are plotted in Figure \ref{medium}. The spatial domain is defined by $\Omega=[0,L_x]\times [0,L_y]$. Specifically, we set  $\Omega=[0,2.2]\times[0,0.6]$ for the first medium and $[0,1]\times[0,1]$ for the second one. The first medium is heterogeneous, i.e., the value changes rapidly in many local regions, which creates great difficulty in approximating the flow dynamics accurately. The second one is a high-contrast medium, where the contrast is approximately 2000. As for the time length, the final time is chosen to be $T=8000$. The fluid properties are set as follows. The relative permeabilities are given by $k_{rw}=\bar{S}_w^2$ and $k_{rw}=(1-\bar{S}_w)^2$, where $\bar{S}_w$ is the effective saturation defined as below:
	\begin{align*}
		\bar{S}_w=\frac{S_w-S_{rw}}{1-S_{rn}-S_{rw}}.
	\end{align*}
	Herein, $S_{r\alpha}$, $\alpha=n,w$ is residual phase saturation and in our experiments and we choose $S_{rw}=S_{rn}=10^{-6}$. Moreover, the viscosity of the wetting phase is set as $\mu_w=1$ and the viscosity for the non-wetting phase is $\mu_n=5$. The density for wetting and non-wetting phases are chosen as $\rho_w=1000$ and  $\rho_n=800$. In this work, we use square fine-scale and coarse-scale elements for spatial discretization. In particular, blocks with side length $h$ are used in the fine-scale partition, where the  mesh size is fixed to be $60\times 220$ for $\kappa_1$ and $100\times 100$ for $\kappa_2$. On the other hand, the coarse mesh size varies, resulting different dimensions of multiscale space. Similarly, there are a fine and coarse partition in time, where the time step size for the former is set to be $\Delta t=100$ and the coarse size may vary. We call solutions computed on the fine mesh in space and smaller time step  reference solutionss, while the one solved in the multiscale space and a coarser time partition is used for approximation. Two types of sources of wetting phase are chosen, denoted by $q_{w,1}$ and $q_{w,2}$. For the non-wetting phase, we set $q_n=0.$ In particular, $q_{w,1}$ and $q_{w,2}$ are defined as follows,

	\begin{align}
		q_{w,1}=\begin{cases}
			0.2, &x\in \eta_1,\\
			-0.2, & x\in \eta_4,\\
			0,  & otherwise.
		\end{cases}
		\quad
		q_{w,2}=\begin{cases}
			0.2, &x\in \eta_i, i=1,\ldots,4\\
			-0.8, & x\in \eta_5,\\
			0,  & otherwise,
		\end{cases}
		\label{source}
	\end{align}
	where $\eta_i$, $i=1,\ldots,4$ are fine-scale elements at the four corners of $\Omega$. In particular, $\eta_1=[0,h]\times [0,h]$, $\eta_2=[L_x-h,L_x]\times [0,h]$, $\eta_3=[0,h]\times [L_y-h,L_y]$, $\eta_4=[L_x-h,L_x]\times[L_y-h,L_y]$. Moreover $\eta_5$
	is the fine-scale element at the center of $\Omega$. More specifically, $\eta_5=[\frac{L_x-h}{2},\frac{L_x+h}{2}]\times [\frac{L_y-h}{2},\frac{L_y+h}{2}]$. We call $q_{w,1}$ a two-point source while $q_{w,2}$ is called five-point source. 

	In this work, we consider the $L^2$ norm to measure the error. In particular, we  define the saturation error as
	\begin{align*}
		e_s=\frac{\|S_{\text{appro}}-S_{\text{ref}}\|_{L^2(\Omega)}}{\|S_{\text{ref}}\|_{L^2(\Omega)}},
	\end{align*} 
	where we use $S_{\text{ref}}$ and $S_{\text{appro}}$ to denote the reference and approximation saturation solutions.
	
\begin{figure}[!htbp] \centering
		\includegraphics[width=5in]{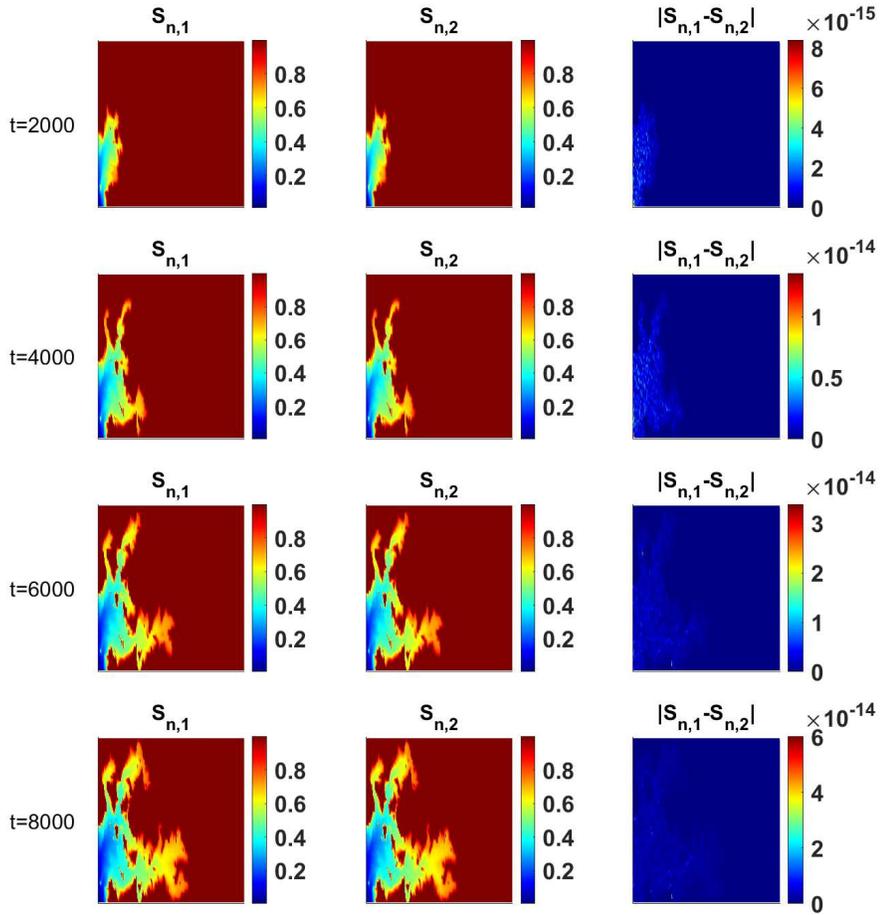}
		\caption{Reference non-wetting phase saturation with $\kappa_1$. We show $S_{n,1}$, $S_{n,2}$ and their absolute differences in three columns. Four time steps $t=2000,4000,6000,8000$ are displayed. Time step size $\Delta t=100$ and fine-scale mesh size is $60\times 220$. }
		\label{conser_ref_k1}
	\end{figure}

\begin{figure}[!htbp] \centering
		\includegraphics[width=5in]{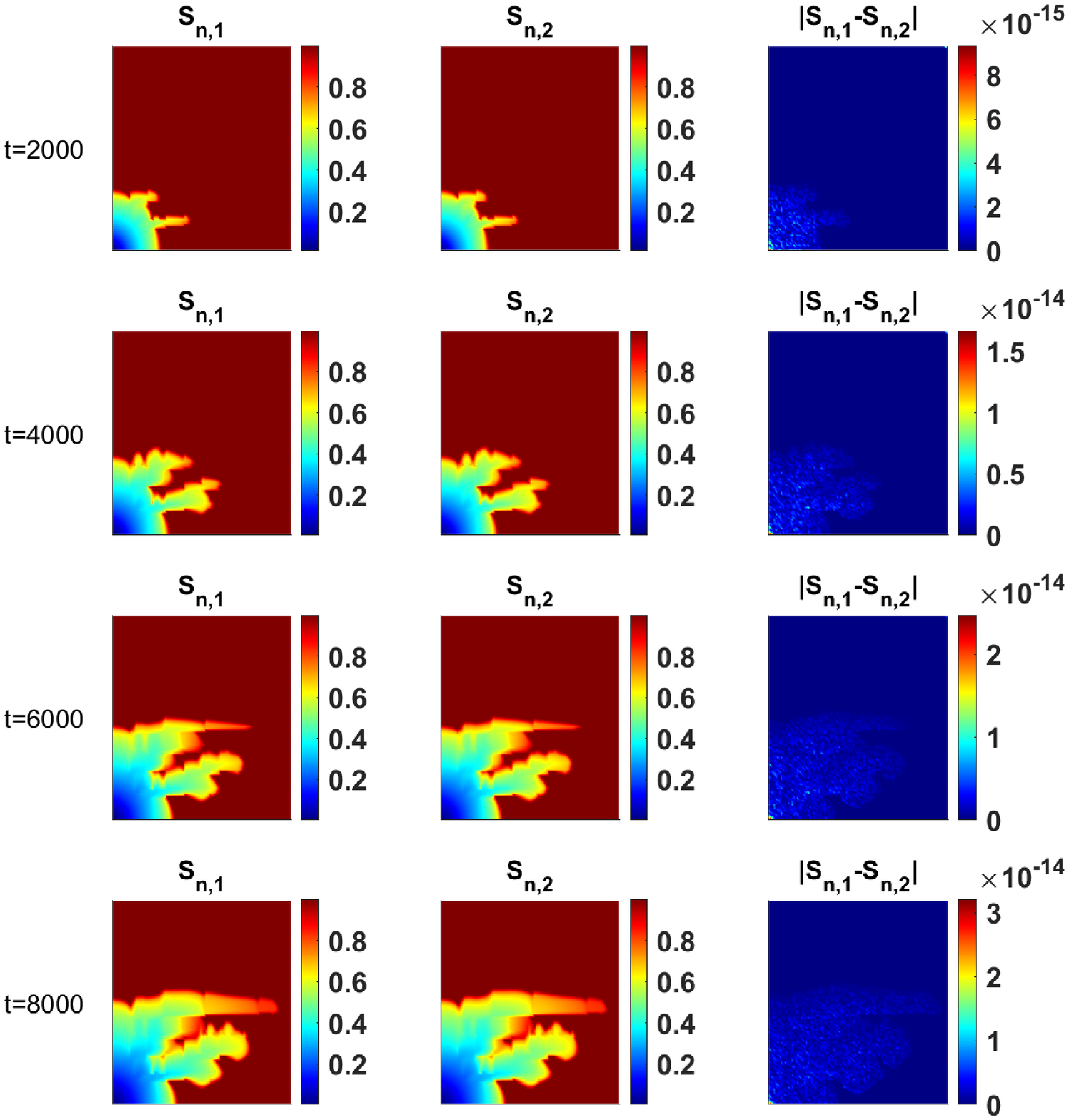}
		\caption{Reference non-wetting phase saturation with $\kappa_2$. We show $S_{n,1}$, $S_{n,2}$ and their absolute differences in three columns. Four time steps $t=2000,4000,6000,8000$ are displayed. Time step size $\Delta t=100$ and fine-scale mesh size is $100\times 100$. }
		\label{conser_ref_k2}
	\end{figure}

\begin{figure}[!htbp] \centering
		\includegraphics[width=5in]{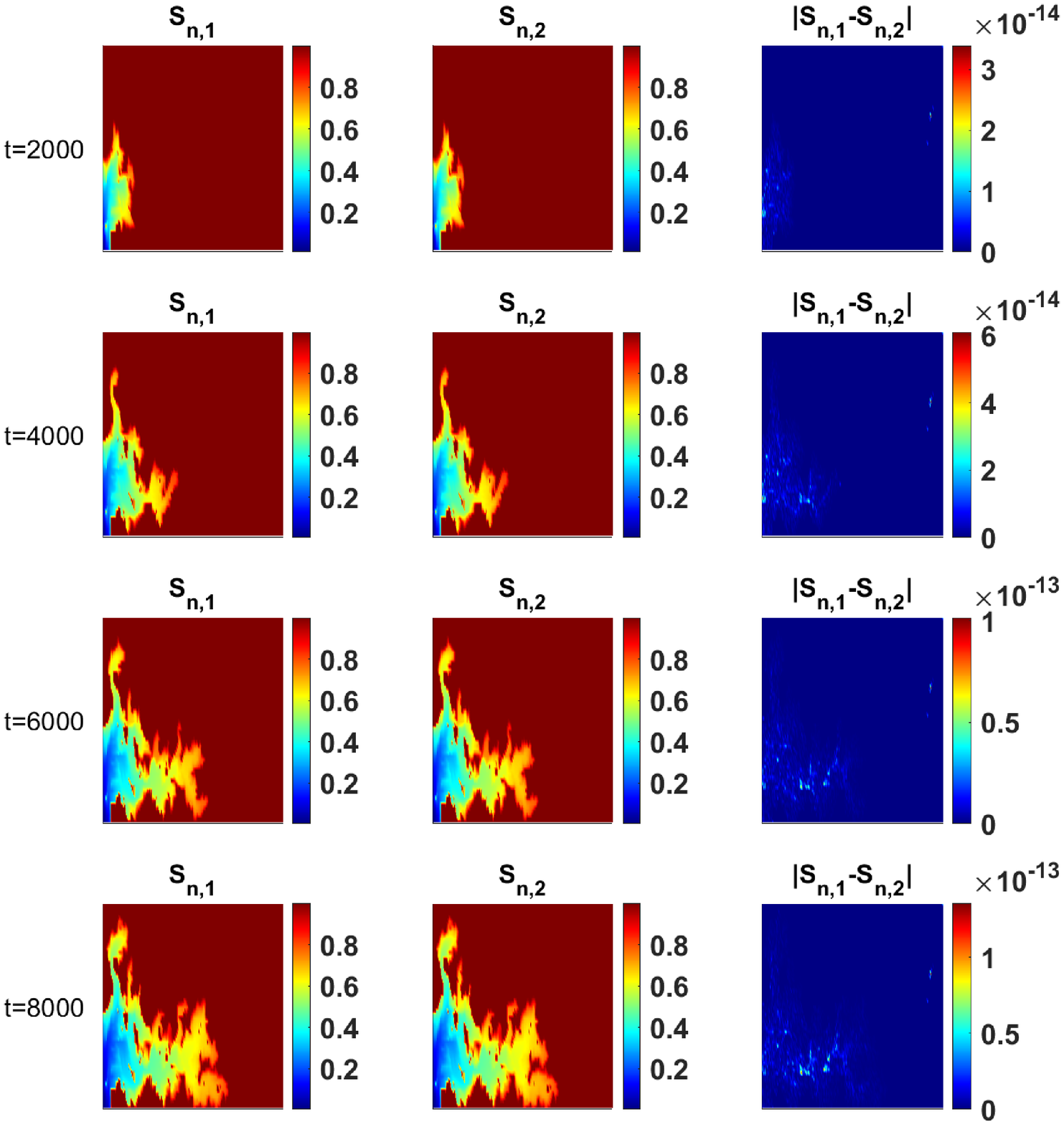}
		\caption{Non-wetting phase saturation computed with ``1+0'' bases with $\kappa_1$. We show $S_{n,1}$, $S_{n,2}$ and their absolute differences in three columns. Four time steps $t=2000,4000,6000,8000$ are displayed. Time step size $\Delta t=100$ and coarse-scale mesh size is $6\times 22$. }
		\label{conser_1+0_k1}
	\end{figure}

\begin{figure}[!htbp] \centering
		\includegraphics[width=5in]{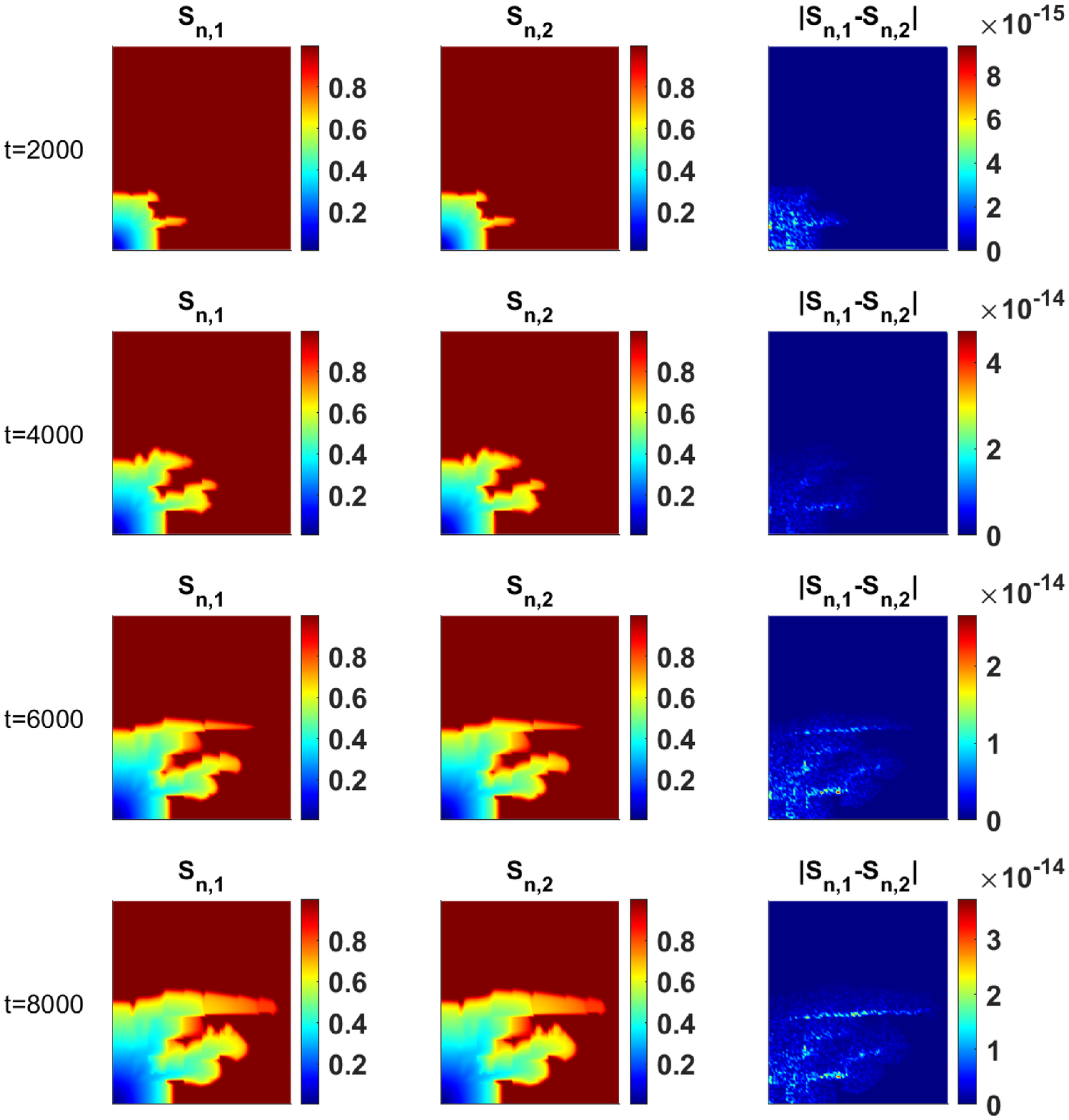}
		\caption{Non-wetting phase saturation computed with ``1+0'' bases with $\kappa_2$. We show $S_{n,1}$, $S_{n,2}$ and their absolute differences in three columns. Four time steps $t=2000,4000,6000,8000$ are displayed. Time step size $\Delta t=100$ and coarse-scale mesh size is $10\times 10$. }
		\label{conser_1+0_k2}
	\end{figure}

	Based on the above notations and settings, we begin with showing the local mass of conservation fulfilled by the results. As mentioned in the subsection \ref{sect:conser}, there are two ways to obtain the non-wetting phase saturation, denoted by $S_{n,1}$ and $S_{n,2}$. If $S_{n,1}=S_{n,2}$, the local mass conservation is attained for both two phases. Hence, we will verify both reference solutions and coarse-scale solutions have fulfilled this requirement by showing that $S_{n,1}$ and $S_{n,2}$ are the same  up to floating-point precision in all the picking time points. In Figures \ref{conser_ref_k1} and \ref{conser_ref_k2}, we show reference non-wetting phase saturation $S_{n,1}$ and $S_{n,2}$ corresponding to $\kappa_1$ and $\kappa_2$. We choose four time steps during the whole simulation process, i.e. $t=2000, 4000,6000, 8000.$  The differences between two non-wetting saturation solutions are almost zero as one can verify from the third column, which reflects that the reference solutions can achieve the conservation of mass. As for the coarse-scale solutions, we choose $``1+0"$ base for testing. If the mass conservation is achieved in this case, it is sufficient to come to conclusion that multiscale solutions can well satisfy this conservation property since solutions computed with more multiscale bases will obtain higher accuracy, which will be tested in the following results. From Figures  \ref{conser_1+0_k1} and  \ref{conser_1+0_k2}, the comparisons corresponding to $\kappa_1$ and $\kappa_2$ are displayed, where the differences are also approximately zero. Thus, multiscale approximations will not lose the local conservation of mass from the above observations. 

\begin{figure}[!htbp] \centering
	\centering
	\subfigure[Reference velocity]
	{ \includegraphics[width=0.45\textwidth]{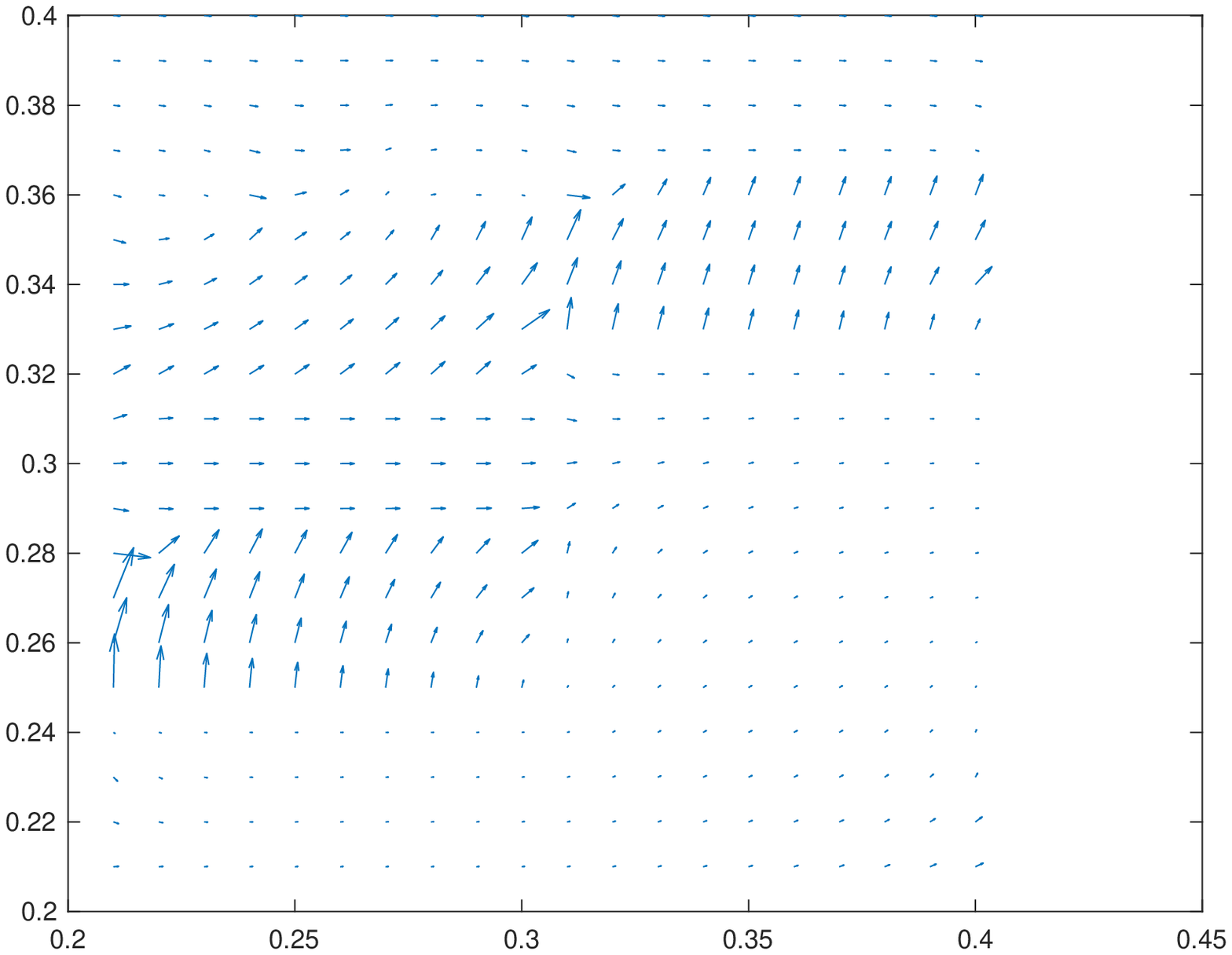}}
	\subfigure[Approximation velocity]
	{ \includegraphics[width=0.45\textwidth]{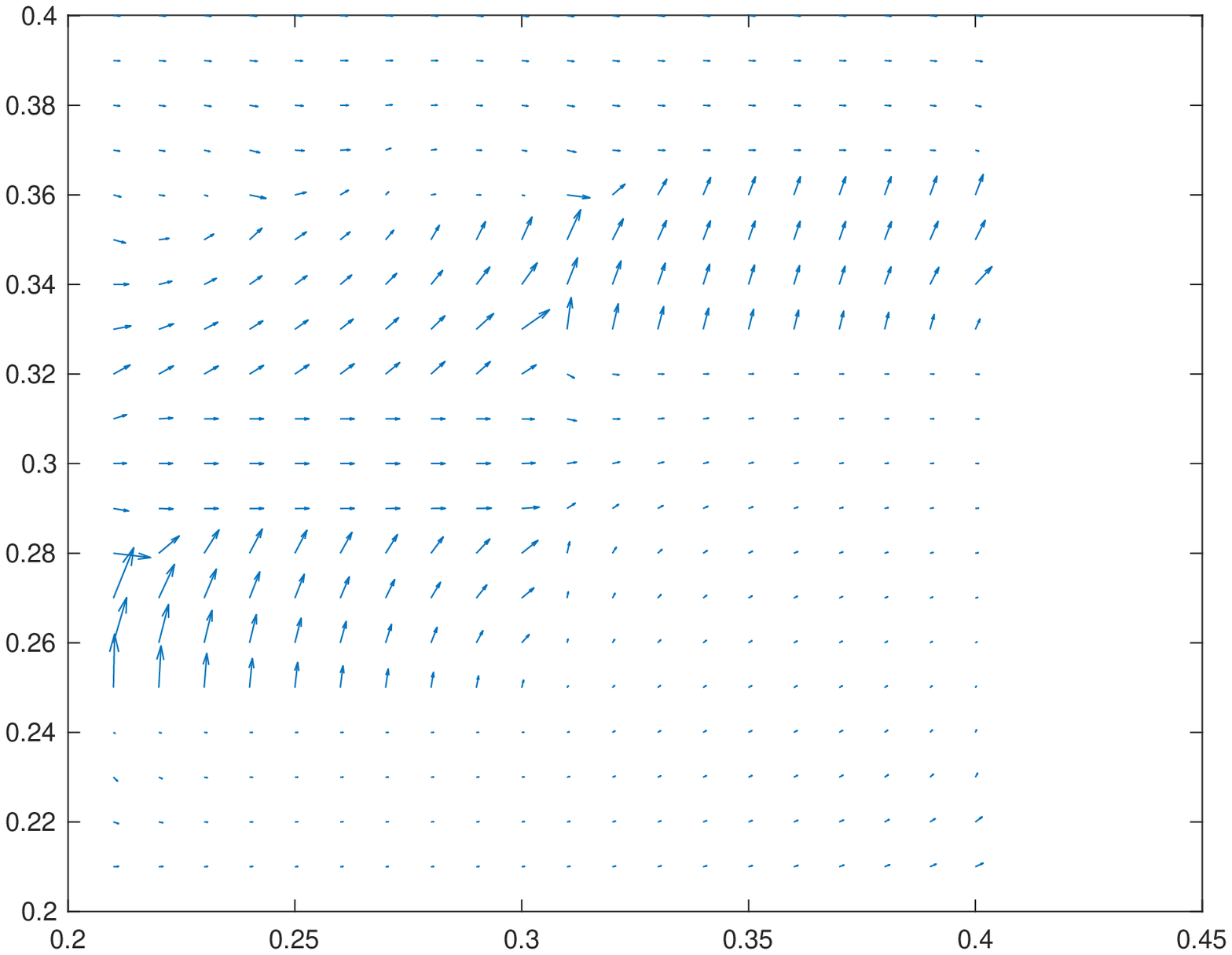}}
	\caption{Comparison of velocity fields with $\kappa_2$. Left: Reference velocity field; Right: approximation velocity field computed with $``2+1"$ multiscale bases. }\label{fig:velo}
\end{figure}

In Figure \ref{fig:velo}, we compare a reference velocity and an approximation field corresponding to the medium $\kappa_2$, where the latter one is computed with $``2+1"$ bases. We remark that a local region is chosen for a clear illustration. Here one can easily observe the approximation is  indistinguishable from the reference. In particular, two solutions share almost the same directions at corresponding places. Hence, Assumption \ref{ass:fw} can be verified.

\begin{figure}[!htbp] \centering
		\subfigure[$\kappa_1$]
		{\includegraphics[width=0.48\textwidth]{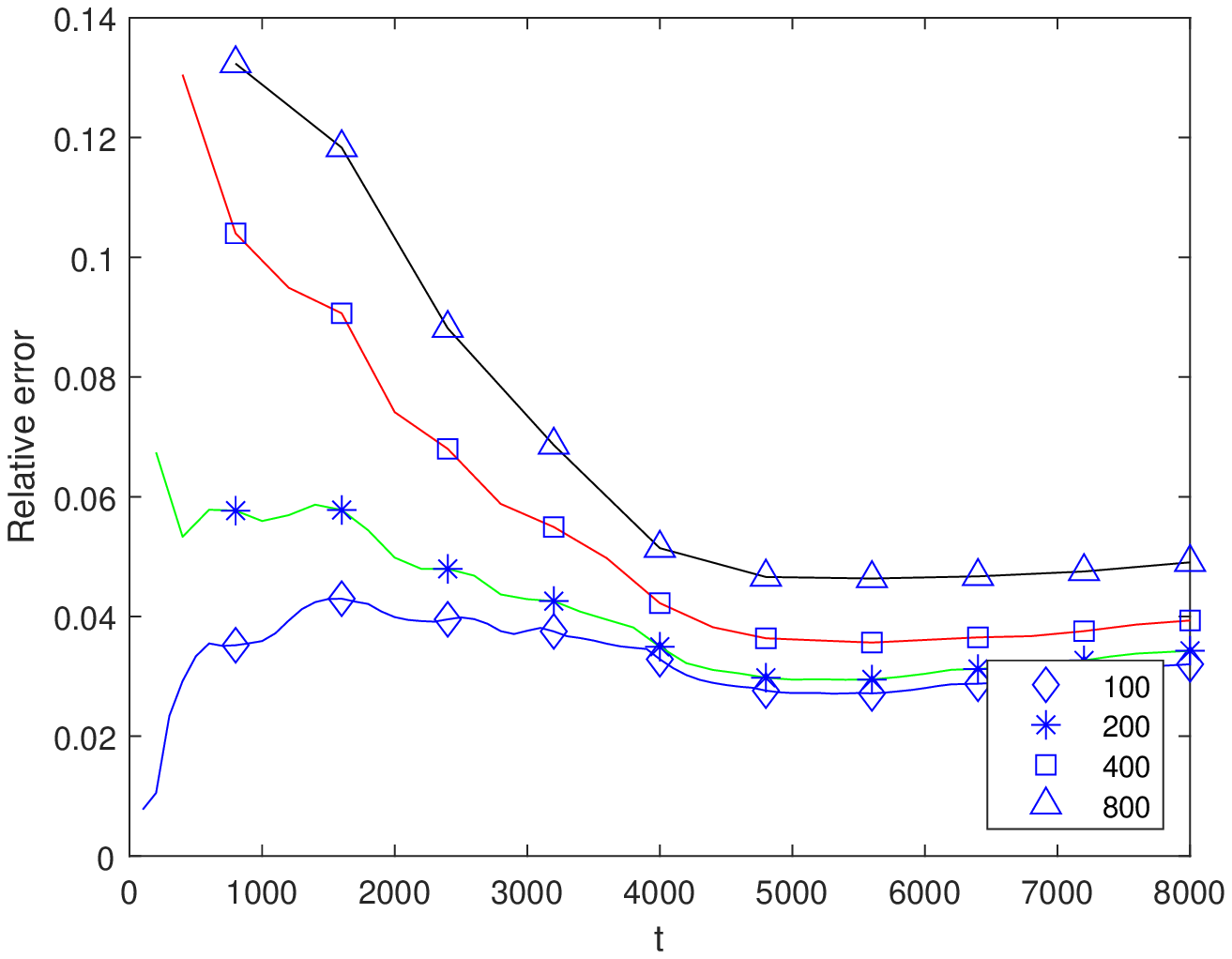}\label{k1_comare_dt}}
		\subfigure[$\kappa_2$]
		{\includegraphics[width=0.48\textwidth]{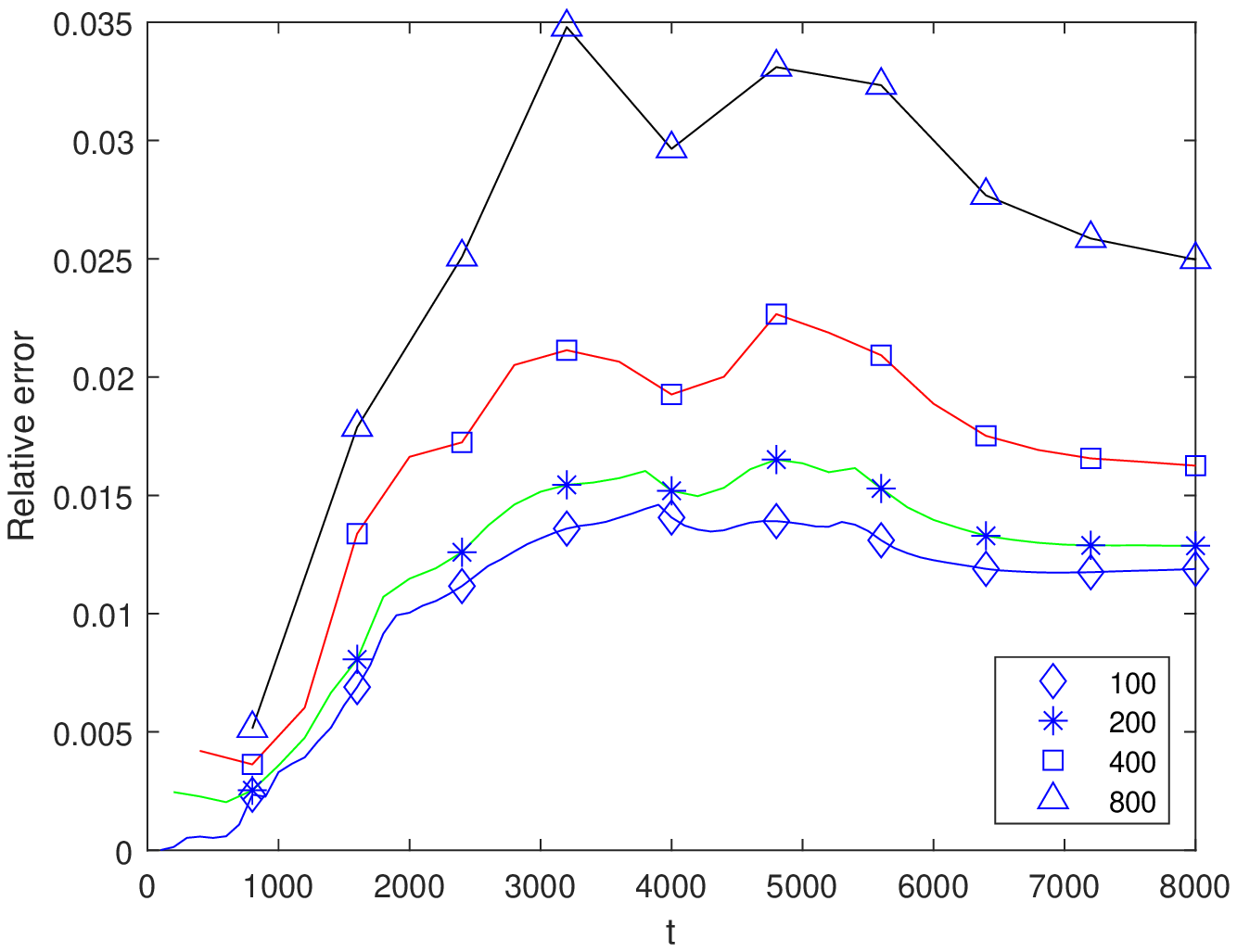}\label{k2_comare_dt}}		
		\caption{Comparison of the $L^2$ errors of the saturation with different time sizes for $\kappa_1$ (left) and $\kappa_2$ (right). We compare four time steps $\Delta t=100,200,400$ and $800$. The coarse mesh sizes are $6\times 22$  for $\kappa_1$ (left) and $10\times 10$ for $\kappa_2$ (right), respectively. The multiscale number is $``5+0"$.}\label{error_dt}
	\end{figure}

    We proceed to demonstrate the error convergence with time step size, mesh size and number of multiscale bases. In the following experiments, it is sufficient to focus on dynamics of wetting-phase saturation owing to the conservation of mass. Thus, without specified, saturation we discuss about refers to that of wetting-phase. We recall that the reference solutions are computed with the fine-scale mesh and with time step $\Delta t=100$ for comparison. In Figure \ref{error_dt}, we present saturation errors with respect to different time step sizes. The subplot on the left corresponds to $\kappa_1$ and the error dynamics of $\kappa_2$ is displayed on the other side. We fix the mesh size and number of multiscale bases in each local region.  The coarse mesh sizes are $6\times 22$ and $10\times 10$ for $\kappa_1$ and $\kappa_2$. In other words, each coarse element is of the same size for $\kappa_1$ and $\kappa_2$, which is composed of $100$ fine-scale elements. We choose $``5+0"$ multiscale bases in each local regions. Four time step sizes are chosen from $100,200,400$ and $800$. Overall, the errors decrease with smaller time step sizes although four curves display different patterns.  With $\kappa_1$ and  $\Delta t=800,400,200$, the errors keep decreasing as time and get steady at about $t=5000$. However, the case with smallest $\Delta t$ is slightly different. The error first increases to $t=2000$  then decreases to $t=5000$ and finally turns fixed. It is worth mentioning that because of updating bases at the intermediate time, i.e. $t=4000$, errors in four curves decease at this time point. With $\kappa_2$, the phenomenon is similar, where the smallest error is obtained with the finest time partition. Besides, one can observe that with a larger $\Delta t$ more time is needed before converging . For example, the error with $\Delta t=800$ has a tendency to decay at the final time while the cases with $\Delta t=100$ and $200$ have been steady from $1000$ time steps before the end. This observation can be explained by the fact that reference solutions are computed with the finest time partition, which means approximations from a smaller time step is easier to converge to corresponding references.

\begin{figure}[!htbp] \centering
		\subfigure[$\kappa_1$]
		{\includegraphics[width=0.48\textwidth]{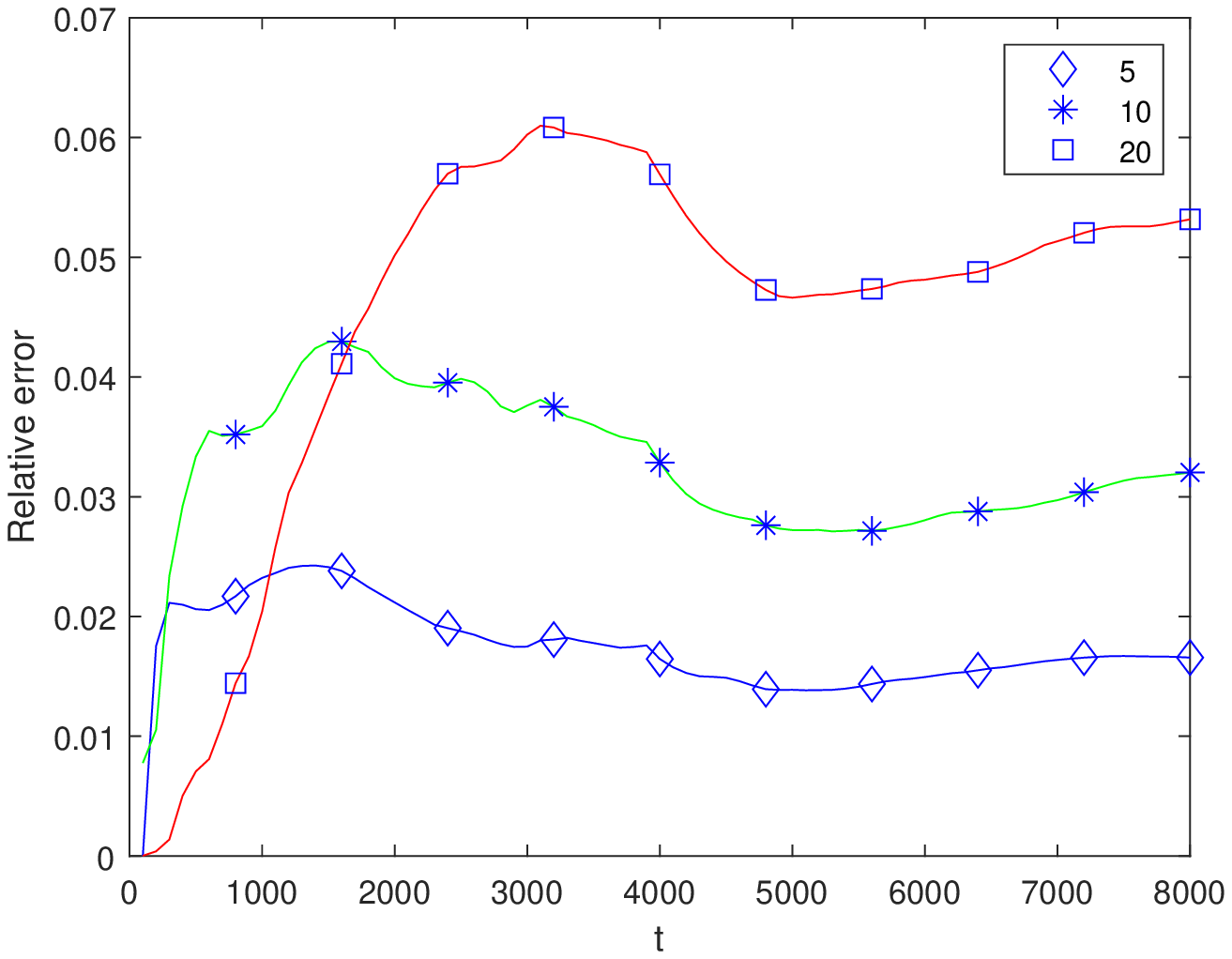}\label{k1_comare_nx}}
		\subfigure[$\kappa_2$]
		{\includegraphics[width=0.48\textwidth]{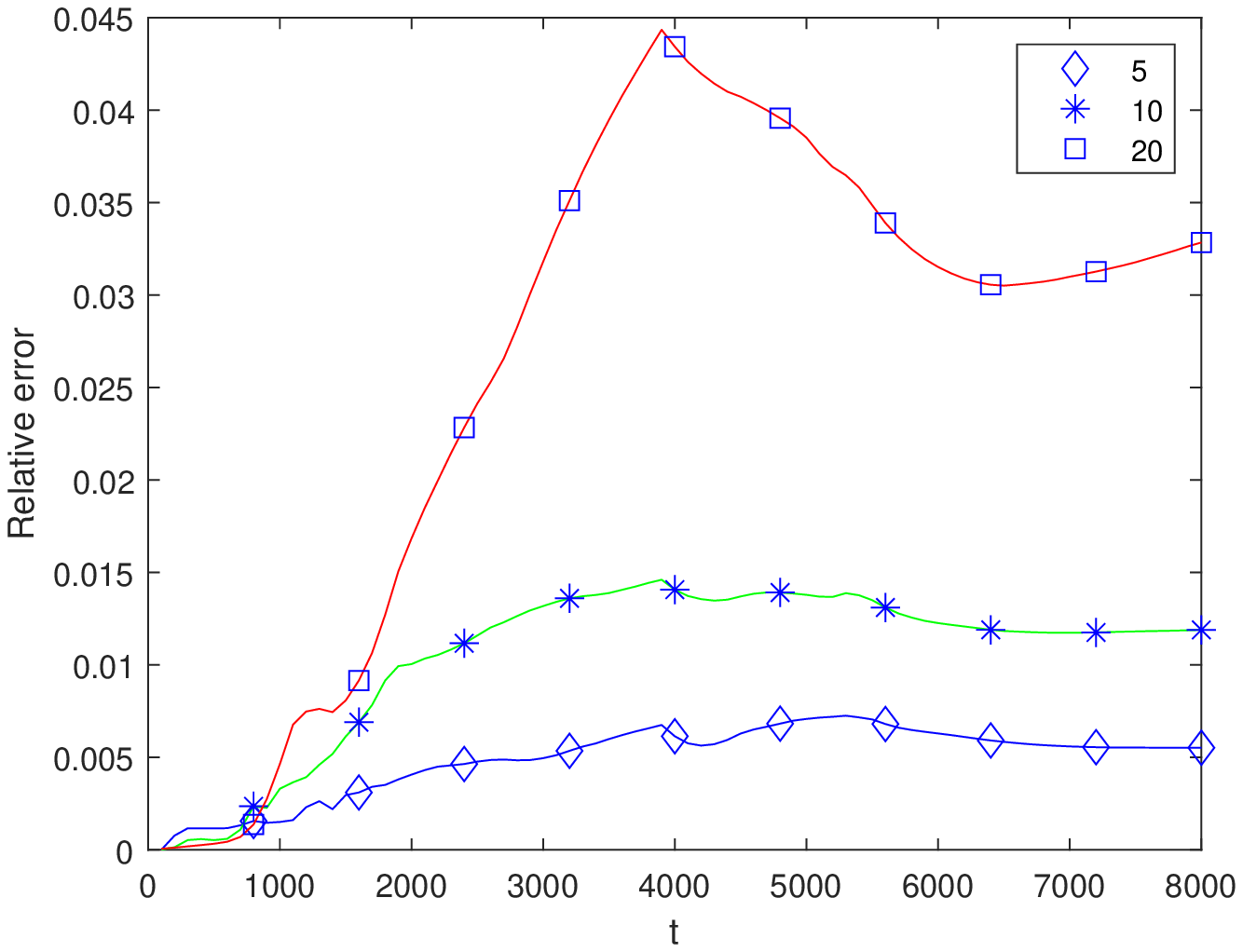}\label{k2_comare_nx}}		
		\caption{Comparison of the $L^2$ errors of the saturation with different mesh sizes for $\kappa_1$ (left) and $\kappa_2$ (right).The mesh size is present by the side length of the coarse square elements, where the legend $n$ means each coarse block is composed of $n\times n$ fine blocks. We compare four mesh sizes $n=5,10$ and $20$. The time step size is $\Delta t=100$ and the multiscale number is $``5+0"$.}\label{error_nx}
	\end{figure}

	In Figure \ref{error_nx}, a comparison is displayed among three coarse mesh sizes. We show the side length of square coarse elements in the legends $n$. For example, $5$ means the blocks contain $5\times5$ small squares in the fine mesh. It is obvious that the approximations are closer to references with smaller coarse elements. From the case with $\kappa_1$, the error is fixed at about $1.5\%$ with smallest coarse blocks while the  error is over $5\%$ at final time when $n=20$. In the right figure, the case is similar, the final error corresponding to $n=20$ (3.2\%)is over six times bigger than that with $n=5$ (0.5\%). Furthermore, it is notable that with smaller coarse elements, the errors converge faster. For instance, in the left figure, the blue curve corresponding to the finest coarse mesh turns smooth at the earliest time. This is similar to the case when we compare different time steps, which reflects that utilizing a finer partition no matter in time or space can speed up the convergence of errors.

\begin{figure}[!htbp] \centering
		\subfigure[$\kappa_1$]
		{\includegraphics[width=0.48\textwidth]{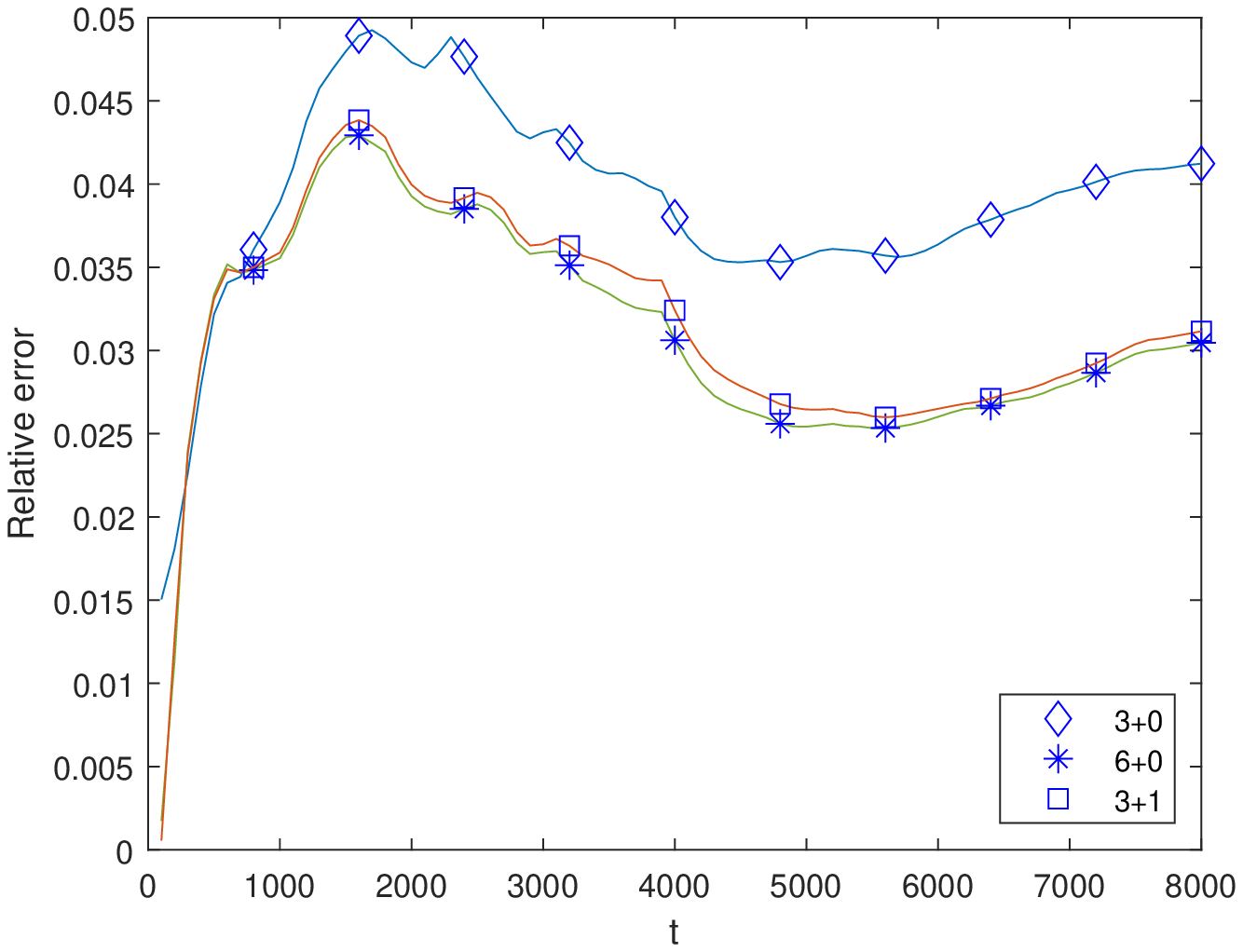}\label{k1_comare_Nb}}
		\subfigure[$\kappa_2$]
		{\includegraphics[width=0.48\textwidth]{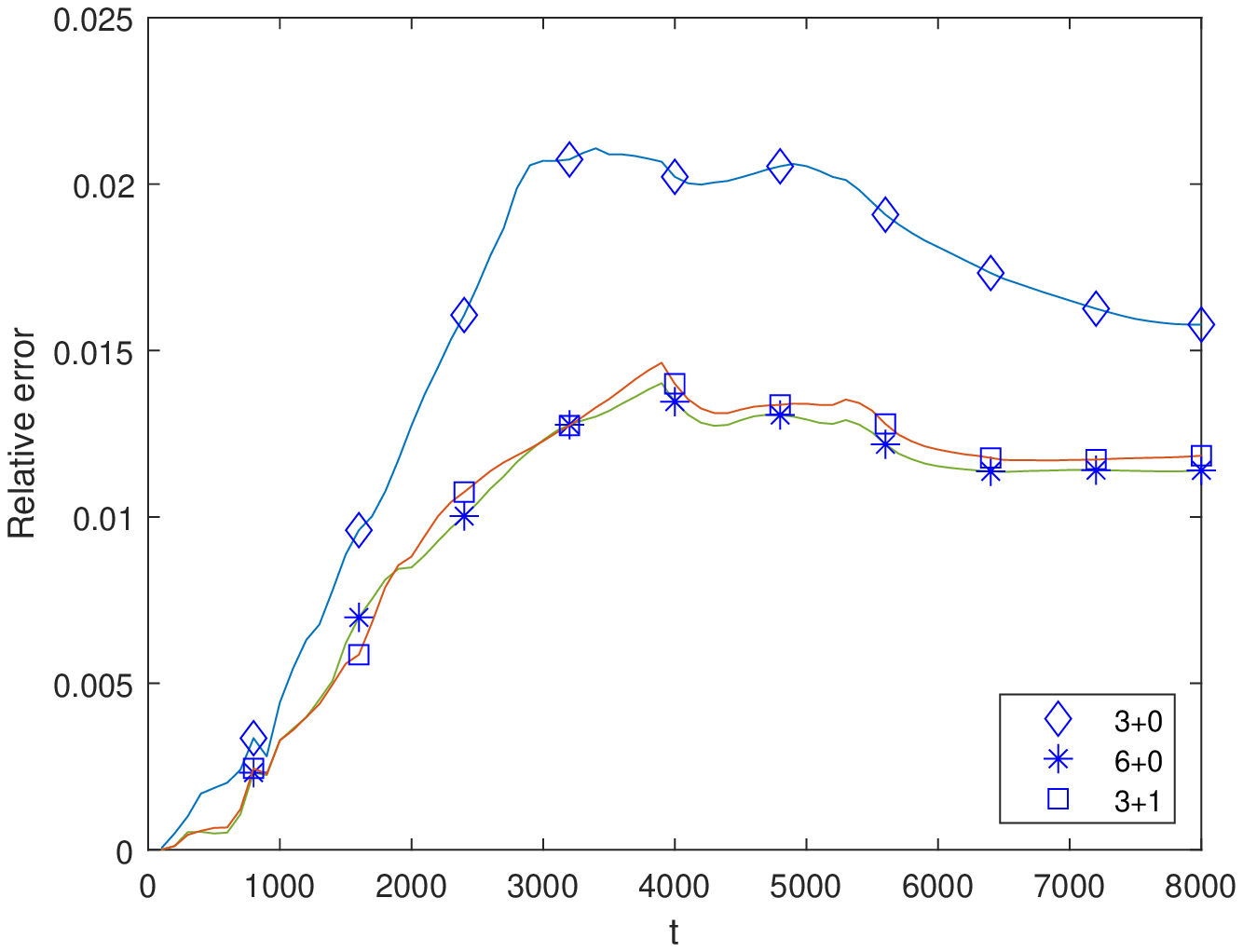}\label{k2_comare_Nb}}		
		\caption{Comparison of the $L^2$ error of the saturation for $\kappa_1$ (left), $\kappa_2$ (right). We compare the approximations with $``3+0"$, $``6+0"$ and $``3+1"$ multiscale bases.  The coarse mesh sizes are $6\times 22$  for $\kappa_1$ (left) and $10\times 10$ for $\kappa_2$ (right), respectively. The time step size is fixed to be $\Delta t=100$. }\label{error_Nb}
	\end{figure}

	Furthermore, we explore the effect of the number of multiscale bases. We fix $\Delta t=100$ and use coarse-scale elements composed of 100 fine-scale elements. Based on this setting, we try three choices, $``3+0"$, $``6+0"$ and $``3+1"$ bases, which are compared in Figure \ref{error_Nb}. Two conclusions can be achieved as follows. First, using more bases in each local region can contribute to lower errors as both $``6+0"$ and $``3+1"$ curves present higher accuracy than the case $``3+0"$. However, the effects of adding offline bases and residual-driven bases are apparently different since adding only one residual-driven basis attains almost the same precision as including three more offline bases in each local region, which shows the efficiency of residual-driven bases in reducing errors. In addition, as we enrich the multiscale bases halfway, one can notice there is a drop in error at this point, which is sharper with $\kappa_1$ than $\kappa_2$. Because the medium $\kappa_1$ shows a higher level of complexity, updating multiscale bases depending on medium may contribute more in increasing accuracy.

\begin{figure}[!htbp] \centering
		\includegraphics[width=5in]{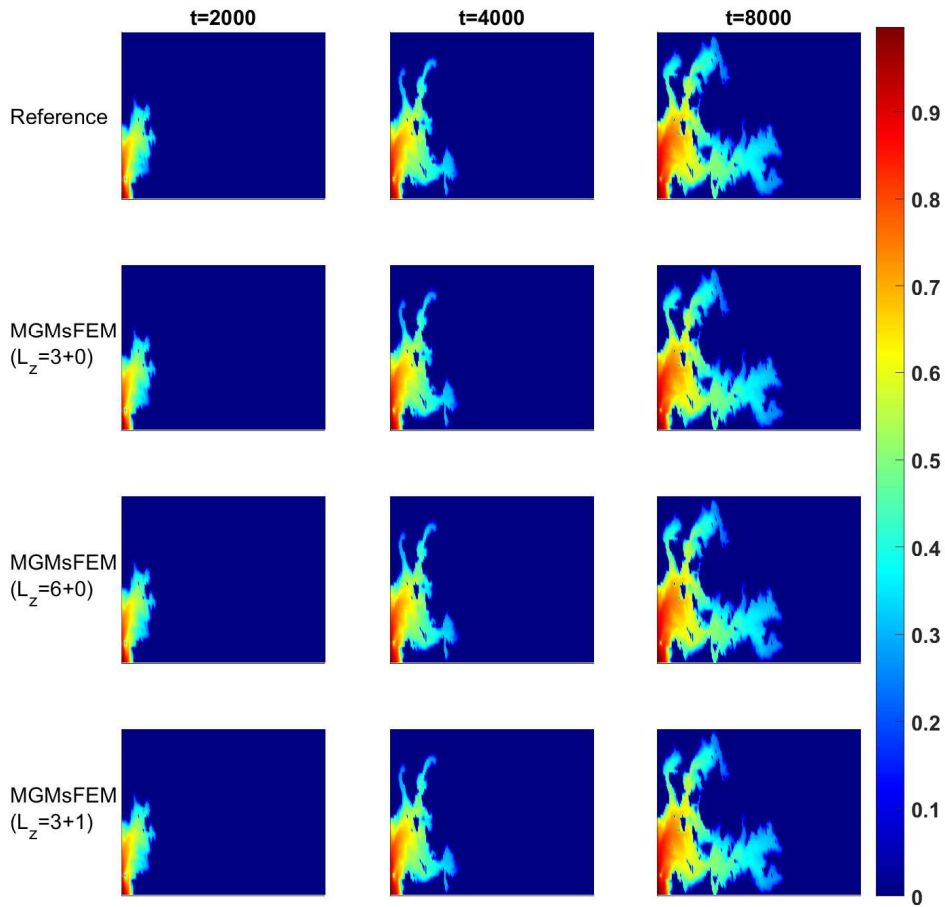}
		\caption{Saturation dynamics with $\kappa_1$. The reference saturation is shown on the top row at three different time levels. The second and third rows are  $L_{ z}=3+0$, $L_z=6+0$, respectively. The last row is using three offline bases and enriched by one online basis for each coarse neighborhood, which is denoted by $L_z=3+1$.}
		\label{s_k1}
	\end{figure}

\begin{figure}[!htbp] \centering
		\includegraphics[width=5in]{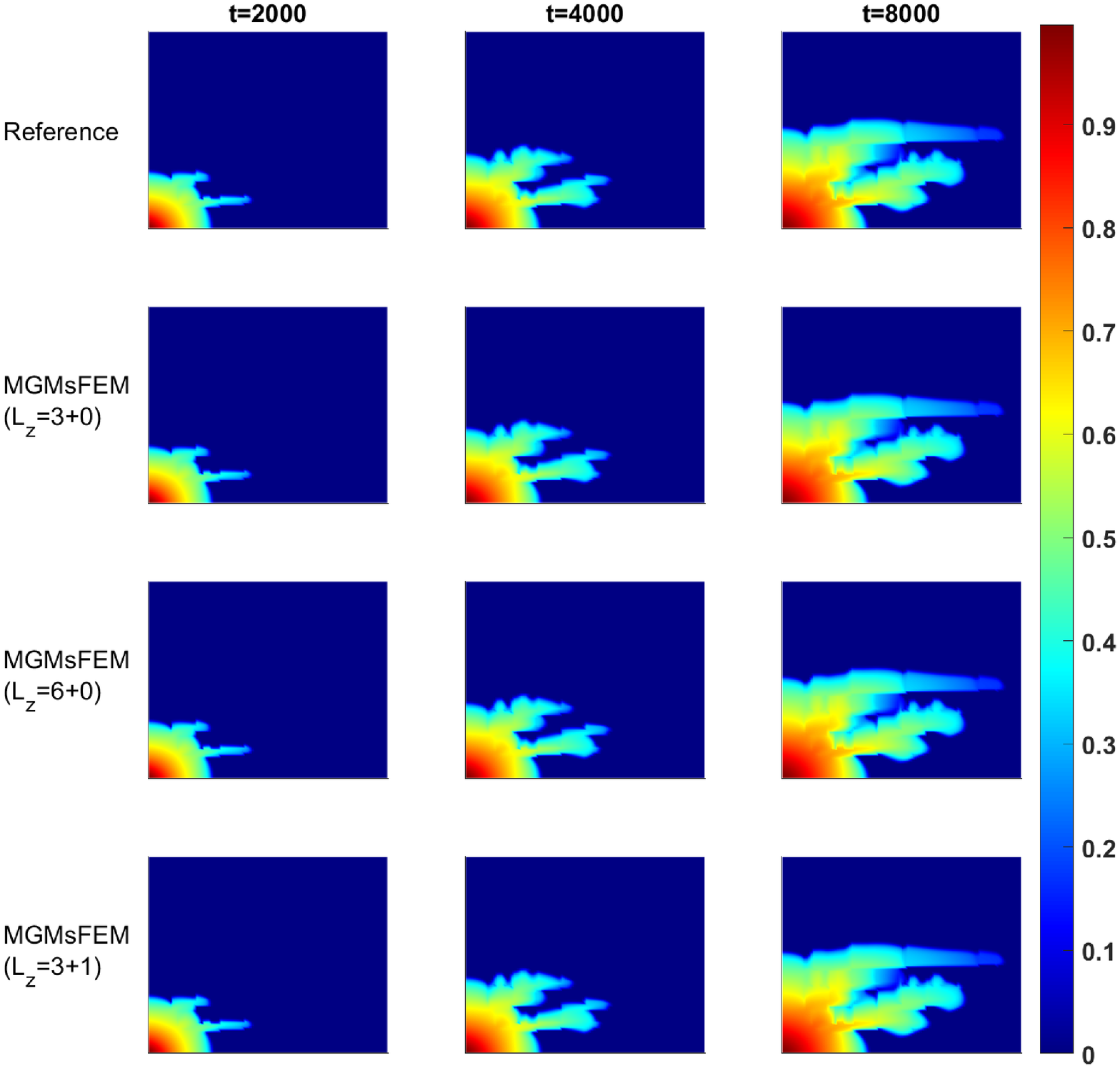}
		\caption{Saturation dynamics with $\kappa_2$. The reference saturation is shown on the top row at three different time levels. The second and third rows are  $L_{ z}=3+0$, $L_z=6+0$, respectively. The last row is using three offline bases and enriched by one online basis for each coarse neighborhood, which is denoted by $L_z=3+1$.}
		\label{s_k2}
	\end{figure}

\begin{figure}[!htbp] \centering
		\includegraphics[width=5in]{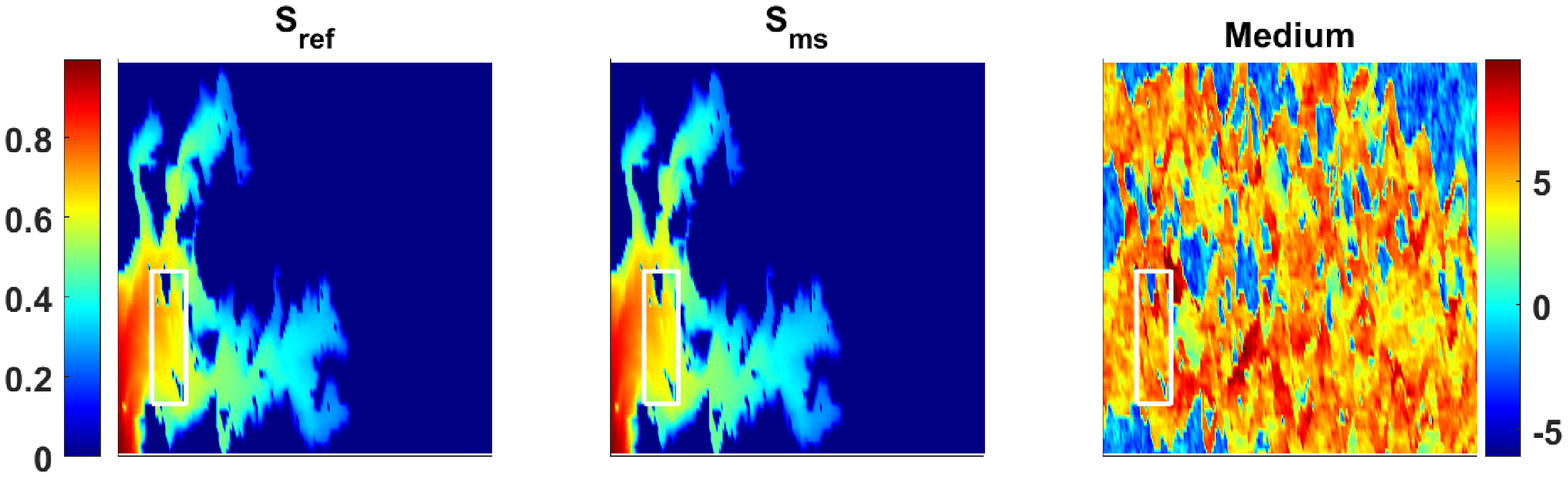}
		\caption{Relationship between medium $\kappa_1$ and corresponding saturation. Left: reference saturation; Middle: approximation saturation; Right: $\kappa_1$. A white block is plotted to highlight the regions where the function changes rapidly. The approximation is solved with $\Delta t=100$ and a coarse mesh size $6\times 22$, where each coarse element is composed of 100 fine blocks. The number of multiscale bases are $``3+1"$. }
		\label{locs_k1}
	\end{figure}

\begin{figure}[!htbp] \centering
		\includegraphics[width=5in]{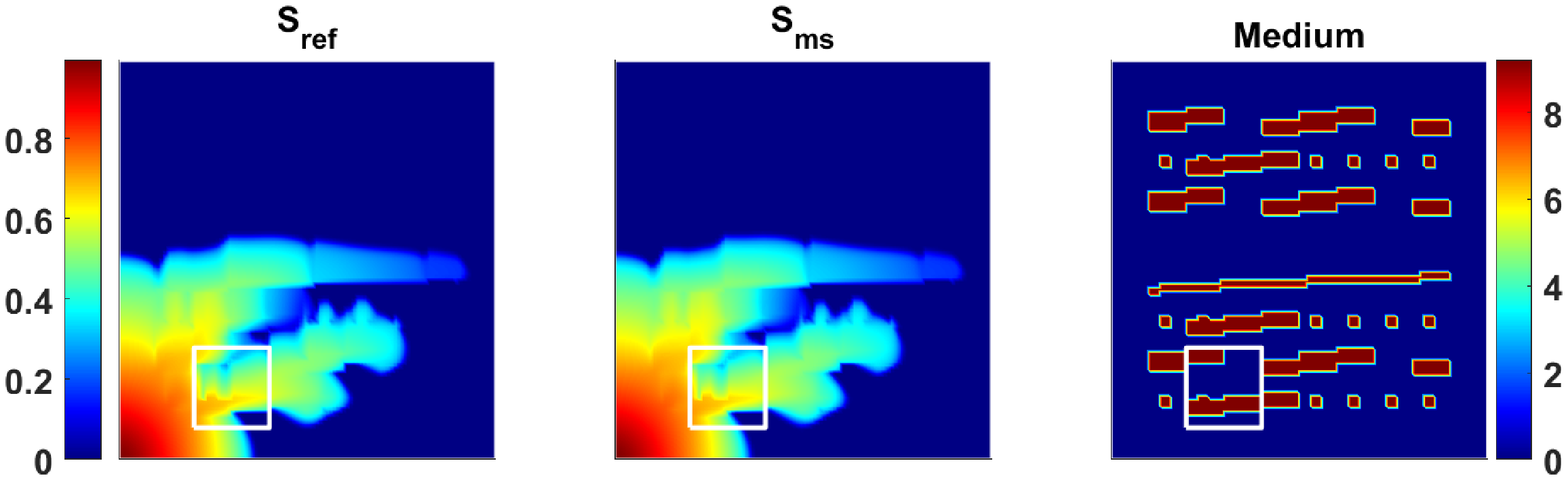}
		\caption{Relationship between medium $\kappa_2$ and corresponding saturation. Left: reference saturation; Middle: approximation saturation; Right: $\kappa_2$. A white block is plotted to highlight the regions where the function changes rapidly. The approximation is solved with $\Delta t=100$ and a coarse mesh size $10\times 10$, where each coarse element is composed of 100 fine blocks. The number of multiscale bases are $``3+1"$. }
		\label{locs_k2}
	\end{figure}

	In Figures \ref{s_k1} and \ref{s_k2}, we show dynamics of reference and approximation saturation at three time points, i.e. $t=2000,4000$ and $8000$. As errors can only indicate accuracy from a global perspective, we can focus on approximation effect of the proposed method in local regions through comparing saturation directly. We demonstrate reference saturation as well as the results from using $``3+0"$, $``6+0"$ and $``3+1"$ multiscale bases. The mesh size and time step are fixed, which are set to be $10\times 10$ coarse elements and $\Delta t=100$. In Figures \ref{s_k1} and \ref{s_k2}, we show the results with respect to $\kappa_1$ and $\kappa_2$. By comparing the approximations with the reference solutions, we can hardly see any obvious distinctions even for the case with highest error (with $``3+0"$ bases). For example, on the first row of Figure \ref{s_k1}, we can see there are blurred regions at the boundary of the reference flow, which are characterized by the solutions computed with multiscale solutions. Hence the fine-scale pattern is well-captured by the approximations, which is shown in Figures \ref{locs_k1} and \ref{locs_k2} in a more straight-forward manner. We still utilize the previous coarse mesh size and time step, i.e. $\Delta t=100$ and coarse elements including $100$ fine blocks.  In both figures, we pick and highlight a local region with a white block where there is an abrupt change in the concerned quantities. Hence, the heterogeneity in medium can result a fast change in the corresponding saturation, which is challenging for approximations. However, the solutions computed with multicsale bases can well characterize the fine-scale pattern, which reflects the high fidelity of the proposed method.

\begin{figure}[htbp!] \centering
		\includegraphics[width=5in]{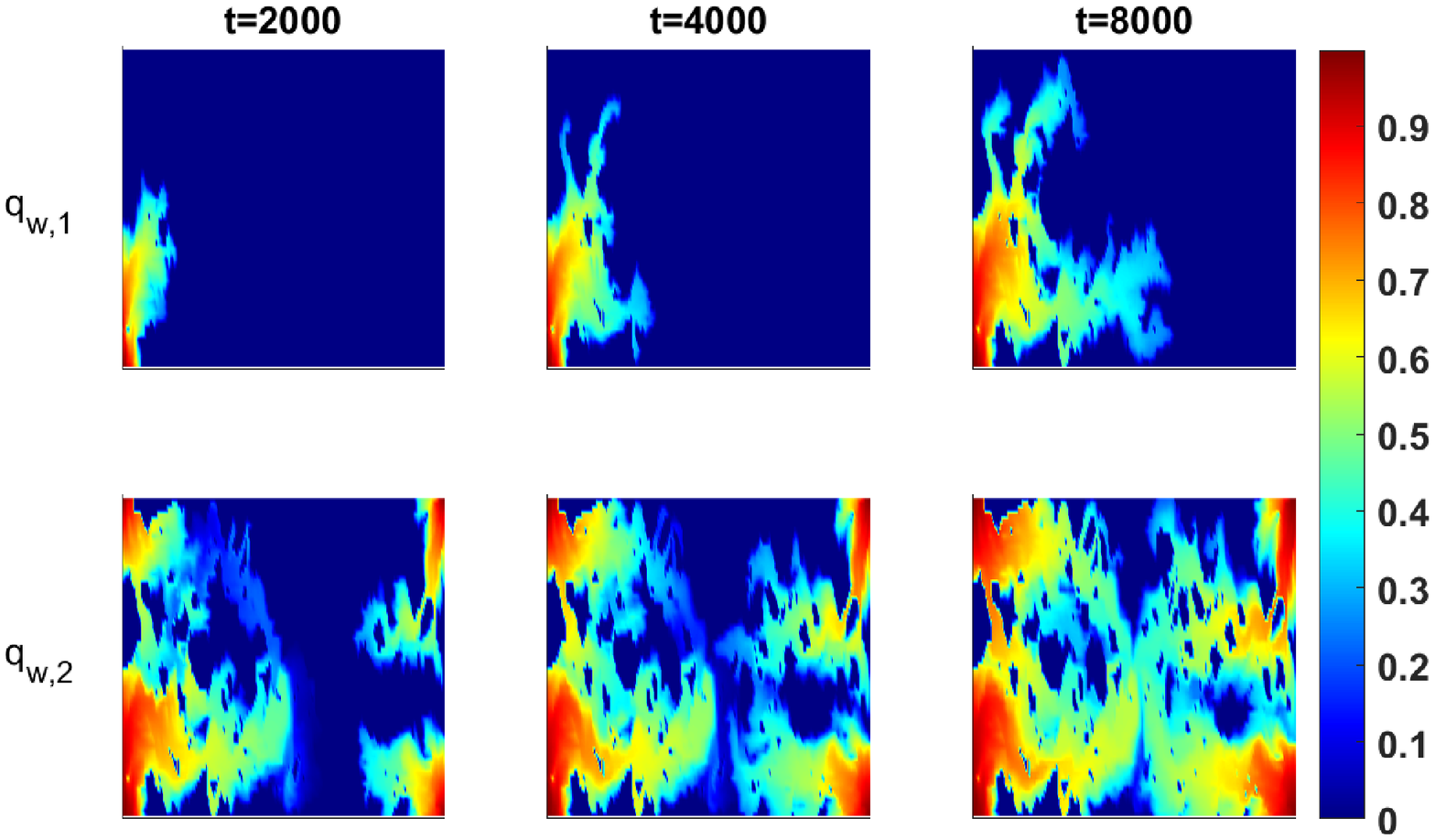}
		\caption{Saturation dynamics with two sources and $\kappa_1$. The time step is fixed to be $\Delta t=100$ and coarse mesh size is $6\times 22$. Reference saturation is used. }
		\label{src_k1}
	\end{figure}

	Finally, we investigate the effect of source terms on saturation dynamics. In Figures \ref{src_k1} and \ref{src_k1}, we plot reference saturation dynamics driven by sources $q_{w,1}$ and $q_{w,2}$ with $\kappa_1$ and $\kappa_2$, respectively. The time step size is set to be $\Delta t=100$ and the coarse mesh size is $6\times 22$ for $\kappa_1$ and $10\times 10$ for $\kappa_2$. Two sources $q_{w,1}$ and $q_{w,2}$ are given in \eqref{source}. We recall that the force values are the same but the locations are different, which creates different effects on the flow motion. In particular, flow moves faster with the five-point source at lower left corner compared with that with the two-point source. Besides, because of the heterogeneous media $\kappa_1$, the dynamics at each corner present different pattern, which is more complicated than the two-point source case.  More specifically, the flow moves faster at two lower corners since the medium values changes more frequently there. On the other hand, the central symmetry in $\kappa_2$ results to a corresponding similar pattern in the flow dynamics, where flows from diagonal corners have similar motions.
	\begin{figure}[htbp!] \centering
		\includegraphics[width=5in]{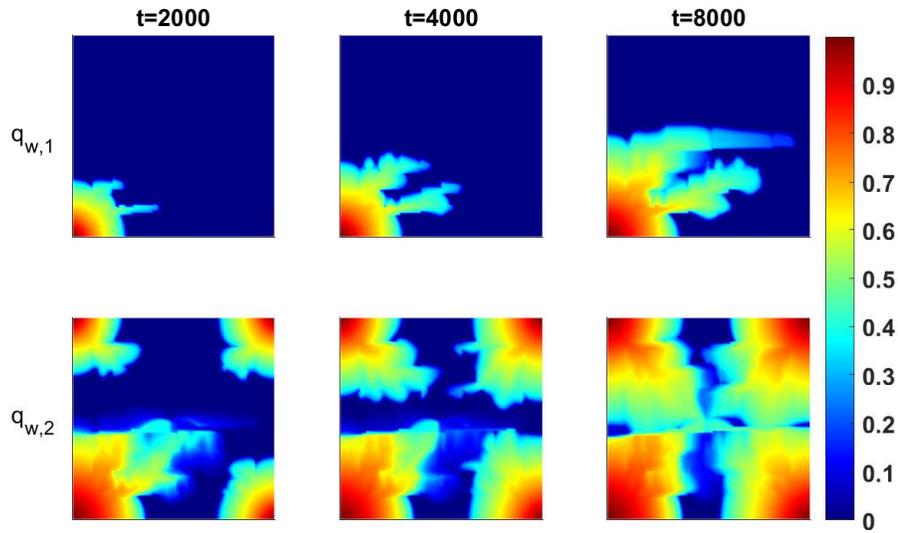}
		\caption{Saturation dynamics with two sources and $\kappa_2$. The time step is fixed to be $\Delta t=100$ and coarse mesh size is $6\times 22$. Reference saturation is used. }
		\label{src_k2}
	\end{figure}
	\section{Analysis}
In this section, we will derive the analysis of the proposed method in two main parts. First, as the analysis in \cite{chen2021new}, we show some important properties of the MS-P-IMPES scheme. Then, we derive  an essential relationship between saturation error and velocity error. In both of analysis and numerical simulations, we neglect the capillary pressure and gravity effect. 

In the first part, we prove three significant properties as follows:
\begin{enumerate}
	\item Local conservation of mass.
	\item Unbiased property of solutions.
	\item bounds-preserving of saturation of both phases.
\end{enumerate}
\subsection{Local conservation of mass and unbiased property}\label{sect:conser}
We will prove that both of the two phases satisfy the local conservation of mass stated as below.
\begin{theorem}
	For any $K^h\in \cT_h$, the saturation of phase $\alpha$ will satisfy the following local conservation of mass as follows,
	\begin{align}
		\int_{K^h} \zeta\frac{S_{\alpha}^{H,n+1}-S_{\alpha}^{H,n}}{\Delta t}+\int_{\pa K^h}f_{\alpha}(S_{w,\alpha}^{*,H,n})\bu_t^{h,n+1}\cdot \bn=\int_{K^h} q_{\alpha}+\sigma_{\alpha}\int_{\pa K^h}f_{w}(S_{w,\alpha}^{*,H,n})f_{n}(S_{w,\alpha}^{*,H,n})\xi_c^{h,n+1}\cdot \bn.
		\label{thm:conser}
	\end{align}
\end{theorem}
\begin{proof}
	To show the first property, we first prove the two schemes to update saturation are equivalent.
	The first scheme is used above in \eqref{ms-scheme2} and \eqref{ms-scheme5}. Here we directly solve saturation of one phase (supposing the wetting phase) and update the other one using that the sum should be equal to $1$.  We review this scheme as follows,
	
	The first saturation update scheme is solving
	\begin{align}
		(\zeta\dfrac{S_{w}^{H,n+1}-S_{w}^{H,n}}{t_{n+1}-t_n},q)+\beta_{w}(\bu_t^{H,n+1},q;S_w^{H,n})&=(q_{w},q)+\beta_c(\xi_{c}^{H,n+1},q;S_w^{H,n}),\label{conser1_Sw}\\
		S_{n}^{H,n+1}&=1-S_{w}^{H,n+1}.	\label{conser1_Sn}	
	\end{align}
	One the other hand, the second scheme solves saturation of both phases directly as below. For all $q\in Q_h$, we solve $S_{\alpha}^{H,n+1}$ for both $\alpha=w,n$ by solving the following two equations. 
	\begin{align}
		(\zeta\dfrac{S_{w}^{H,n+1}-S_{w}^{H,n}}{t_{n+1}-t_n},q)+\beta_{w}(\bu_t^{H,n+1},q;S_w^{H,n})&=(q_{w},q)+\beta_c(\xi_{c}^{H,n+1},q;S_w^{H,n}),\label{conser2_Sw}\\
		(\zeta\dfrac{S_{n}^{H,n+1}-S_{n}^{H,n}}{t_{n+1}-t_n},q)+\beta_{n}(\bu_t^{H,n+1},q;S_w^{H,n})&=(q_{n},q)-\beta_c(\xi_{c}^{H,n+1},q;S_w^{H,n}).\label{conser2_Sn}	
	\end{align}
	Since \eqref{conser1_Sw} is the same as \eqref{conser2_Sw}, it remains to show that $S_{n}^{H,n+1}$ solved by \eqref{conser1_Sn} and \eqref{conser2_Sn} are the same. Due to the postpocessing technique, 	\eqref{ms-scheme1} can be replaced by \eqref{pp}. Subtract \eqref{conser1_Sw} by \eqref{pp} and we can obtain 
	\begin{align*}
		(\zeta\dfrac{S_{w}^{H,n+1}-S_{w}^{H,n}}{t_{n+1}-t_n},q)-\beta_{n}(\bu_t^{H,n+1},q;S_w^{H,n})&=-(q_{n},q)+\beta_c(\xi_{c}^{H,n+1},q;S_w^{H,n}). 
	\end{align*}
	Apply \eqref{conser1_Sn} and we can obtain \eqref{conser2_Sn}.
	Reversely, one can attain \eqref{conser1_Sn} by combining  \eqref{conser2_Sw} and \eqref{conser2_Sn}. Indeed, by adding \eqref{conser2_Sw} to \eqref{conser2_Sn}, we can obtain 
	\begin{align*}
		\sum\limits_{\alpha}(\zeta\dfrac{S_{\alpha}^{H,n+1}-S_{\alpha}^{H,n}}{t_{n+1}-t_n},q)+\sum\limits_{\alpha}\beta_{\alpha}(\bu_t^{H,n+1},q;S_w^{H,n})=(q_{t},q).
	\end{align*}
	Using \eqref{pp}, it holds that
	\begin{align*}
		\sum\limits_{\alpha}(\dfrac{S_{\alpha}^{H,n+1}-S_{\alpha}^{H,n}}{t_{n+1}-t_n},q)=0.
	\end{align*}
	Taking $q=1_{K^h}$ for each $K^h$, we can have 
	\begin{align*}
		\sum\limits_{\alpha}\dfrac{S_{\alpha}^{H,n+1}-S_{\alpha}^{H,n}}{t_{n+1}-t_n}=0
	\end{align*}
	on each $K^h$. Since $S_w^{H,0}+S_w^{H,0}=1$, \eqref{conser1_Sn} can be attained for all $n$. 
	Hence, both schemes are equivalent. Finally, we can obtain \eqref{thm:conser} by choosing $q=1_{K^h}$ in \eqref{conser2_Sw} and \eqref{conser2_Sn}.	
\end{proof}

Then, we show the scheme MS-P-IMPES is unbiased in the solutions $S_{\alpha}^{H,n+1}$, $\bu_{\alpha}^{H,n+1}$ and $p_{\alpha}^{H,n+1}$ for both $\alpha=w,n$.
From both two schemes of updating saturation, we can easily verify that $S_{w}^{H,n+1}$ and $S_{n}^{H,n+1}$ are unbiased. Besides, $\bu_t^{H,n+1}$ and $p_w^{H,n+1}$ are uniquely determined by \eqref{ms-scheme1} and \eqref{ms-scheme3}. Also, $\xi_c^{H,n+1}$ is uniquely solved by \eqref{ms-scheme4}. Hence, $p_{n}^{H,n+1}$ can be obtained by \eqref{ms-scheme6}. Now we remain to show $\bu_{\alpha}^{H,n+1}$ is unbiased for $\alpha=w,n$. It can be verified using the following relations:
\begin{align*}
	\bu_w^{H,n+1}&=f_w(S_{w}^{H,n})\bu_t^{H,n+1}-f_w(S_{w}^{H,n})f_n(S_{w}^{H,n})\xi_c^{H,n+1},\\
	\bu_n^{H,n+1}&=f_n(S_{w}^{H,n})\bu_t^{H,n+1}+f_w(S_{w}^{H,n})f_n(S_{w}^{H,n})\xi_c^{H,n+1}.
\end{align*}
Hence, we finish the proof in the unbiased property of the MS-P-IMPES scheme.

\subsection{bounds-preserving property}
This part is motivated by \cite{chen2021new}. First, we need two lemmas.
\begin{lemma}\label{ass1}
	There exits a positive constant $l_{\alpha}$ such that
	\begin{align*}
		f_{\alpha}(S_{\alpha})\leq l_{\alpha} S_{\alpha},\quad \alpha=w,n.
	\end{align*}
\end{lemma}
This lemma informs that the fractional flow functions are bounded by linear functions of saturation.
\begin{lemma}\label{ass2}
	We define $q_{\alpha}$ in \eqref{model} a sink term if $q_{\alpha}\leq 0$. In this case, there exist two positive constants $l_1$ and $l_2$ such that
	\begin{align*}
		l_1S_{\alpha}\leq |q_{\alpha}|\leq l_2 S_{\alpha}.
	\end{align*}
\end{lemma}
The above two lemmas are proved in \cite{chen2021new}. Based on them, the theorem that guarantees the bounds-preserving property is stated below.
\begin{theorem}
	We assume $S_{\alpha}^{H,n}\in (0,1)$ with tolerance saturation $S_{t\alpha}$ satisfying $S_{\alpha}^{H,n}\geq S_{t\alpha}>0$, $\alpha=w,n$. We define $\Delta t=t_{n+1}-t_n$. For any $K^h\in \cT_h$, if $\frac{\Delta t}{h}$ is small enough, it holds that 
	\begin{align*}
		S_{\alpha}^{H,n+1}\in(0,1), \quad \alpha=w,n.
	\end{align*}
\end{theorem}
\begin{proof}
	The proof is motivated by Lemma 4.3 in \cite{chen2021new} except that the $\bu_t^{h,n+1}$ in the condition of determine proper $\Delta t$ should be replaced by $\bu_t^{H,n+1}$, which is solved in the MS-P-IMPES scheme. We present the proof below for the convenience of readers.
	
	For all $K^h\in \cT_h$, we choose $q=1_{K^h}$ in \eqref{conser1_Sw} and we have
	\begin{align}
		\zeta(S_w^{H,n+1}-S_w^{H,n})=\Delta t q_w+\frac{\Delta t}{h}\sum_{F\subset \pa K^h}f_w(S_{w,w}^{*,H,n})(f_n(S_{w,n}^{*,H,n})\xi_c^{H,n+1}\cdot\bn-\bu_t^{H,n+1}\cdot\bn)|_F. \label{temp1}
	\end{align}
	We define $\bv_w^{H,n+1}\in V_{h}$ by $\bv_w^{H,n+1}\cdot \bn|_F=\left.\left(\bu_t^{H,n+1}\cdot\bn-f_n(S_{w,n}^{*,H,n})\xi_c^{H,n+1}\cdot \bn\right)\right|_F$ on any $F\in \cE_f$. 
	We further define $\pa K^{h,+}_{w,n}=\{F\subset \pa K^h: \bv_w^{H,n}\cdot \bn|_F\geq 0\}$ and $\pa K^{h,-}_{w,n}=\{F\subset \pa K^h: \bv_w^{H,n}\cdot \bn|_F<0\}$.
	Hence  \eqref{temp1} is equivalent to 
	\begin{align*}
		S_w^{H,n+1}&= S_w^{H,n}+\frac{\Delta t }{\zeta}q_w-\frac{\Delta t}{\zeta h}\left(\sum_{F\subset  \pa K^{h,+}_{w,n+1}}+\sum_{F\subset  \pa K^{h,-}_{w,n+1}}\right)f_w(S_{w,w}^{*,H,n})\bv_w^{H,n+1}\cdot \bn|_F.
	\end{align*}
	To prove $S_w^{H,n+1}\geq 0$, we estimate $\sum_{F\subset  \pa K^{h,+}_{w,n+1}}f_w(S_{w,w}^{*,H,n})\bv_w^{H,n+1}\cdot \bn|_F$. We can determine $f_w(S_{w,w}^{*,H,n})$ by partitioning $\pa K^{h,+}_{w,n+1}$ into two groups. On $F\subset  \pa K^{h,+}_{w,n+1}\cap \pa K^{h,+}_{w,n}$, $S_{w,w}^{*,H,n}=S_w^{H,n}|_{K^h}$ and on $F\subset  \pa K^{h,+}_{w,n+1}\setminus \pa K^{h,+}_{w,n}$, $S_{w,w}^{*,H,n}=S_w^{H,n}|_{ K^{h,F}}$ where $K^{h,F}$ is the neighboring element satisfying $K^h\cap K^{h,F}=F$. Then we use Lemma \ref{ass1} and it follows that
	\begin{align*}
		&\sum_{F\subset  \pa K^{h,+}_{w,n+1}}f_w(S_{w,w}^{*,H,n})\bv_w^{H,n+1}\cdot \bn|_F\\
		&=\left(\sum_{F\subset  \pa K^{h,+}_{w,n+1}\cap \pa K^{h,+}_{w,n}}+\sum_{F\subset  \pa K^{h,+}_{w,n+1}\setminus \pa K^{h,+}_{w,n}}\right)f_w(S_{w,w}^{*,H,n})\bv_w^{H,n+1}\cdot \bn|_F\\
		&\leq \sum_{F\subset  \pa K^{h,+}_{w,n+1}\cap \pa K^{h,+}_{w,n}}l_wS_w^{H,n}|_{K^h}\bv_w^{H,n+1}\cdot \bn|_F+
		\sum_{F\subset  \pa K^{h,+}_{w,n+1}\setminus \pa K^{h,+}_{w,n}}l_wS_w^{H,n}|_{K^{h,1}}\bv_w^{H,n+1}\cdot \bn|_F\\
		&\leq \sum_{F\subset  \pa K^{h,+}_{w,n+1}}l_wS_w^{H,n}|_{K^h}\bv_w^{H,n+1}\cdot \bn|_F+
		\sum_{F\subset  \pa K^{h,+}_{w,n+1}\setminus \pa K^{h,+}_{w,n}}l_w\left|S_w^{H,n}|_{K^{h,1}}-S_w^{H,n}|_{K^{h}}\right|\bv_w^{H,n+1}\cdot \bn|_F.
	\end{align*}
	We recall that $S_{w}^{H,n}\geq S_{tw}>0$, hence there exists  $\epsilon_1$ and $\epsilon_2$ such that when $\frac{\Delta t}{\zeta h}\leq \epsilon_1$ and  $|\bv_w^{H,n+1}\cdot \bn|_F\leq \epsilon_2$ with $F\subset  \pa K^{h,+}_{w,n+1}\setminus \pa K^{h,+}_{w,n}$, it holds that $S_{w}^{H,n+1}\geq \eta S_{tw}>0$ for some constant $\eta\in (0,1)$.  In particular, $\Delta t$ should be small enough such that the following two requirements satisfy:
	\begin{enumerate}
		\item $\frac{\Delta t}{\zeta h}\leq \epsilon_1$;
		\item $|\bv_w^{H,n+1}\cdot \bn|_F\leq \epsilon_2$ with $F\subset  \pa K^{h,+}_{w,n+1}\setminus \pa K^{h,+}_{w,n}$.
	\end{enumerate}
	We remark that the second condition requires $|\bv_w^{H,n+1}\cdot \bn|$ to be small enough on the fine edges where $\bv_w^{H,n+1}\cdot \bn$ changes sign from $t_n$ to $t_{n+1}$, which can be achieved if the $\Delta t$ is small enough. We further define $\Delta t_1$ to be the smallest value satisfying the above two requirements. Hence with $\Delta t \leq \Delta t_1$, we have $S_{w}^{H,n+1}\geq \eta S_{tw}>0$.  Similarly, we can determine $\Delta t_2$ such that when $\Delta t\leq \Delta t_2$, $S_{n}^{H,n+1}\geq \eta S_{tn}>0$. Due to the relation $S_{w}^{H,n+1}+S_{n}^{H,n+1}=1$, we can easily obtain that $S_{\alpha}^{H,n+1}<1$ for $\alpha=w,n.$ Finally, we  choose $\Delta t=\min\{\Delta t_1,\Delta t_2\}$ and then $S_{\alpha}^{H,n+1}\in(0,1)$.
\end{proof}

\subsection{Error analysis}
In this part, we derive the error analysis for saturation. For simplification, we consider only consider linear or quadratic relative permeability, i.e. $k_{r\alpha}=S_{\alpha}^p$, for $p=1,2$. Without loss of generality, only in this part, we assume $\mu_w=\mu_n=1$ since we can always put the constant viscosity into the constant coefficients of the final errors. Besides, we choose $p_{\alpha}^B=0$ and $g_{\alpha}^N=0$ for $\alpha=w,n$. We also assume that $\kappa_{\text{min},0}=\min_{\bx} K(\bx)=1$. In this way, our proof will be present in a more concise way.

\begin{assumption}\label{ass:fw}
For each $K^h$, we assume that on each $e\subset \pa K^h$, it holds that
\begin{align*}
	(u_t^{h,n+1}\cdot \bn)(u_t^{H,n+1}\cdot \bn)|_e\geq 0.
\end{align*}
\end{assumption}
This assumption states that the reference  and approximation flows should have the same direction on each fine edge and it can be verified by Figure \ref{fig:velo}.
We then proceed to the relationship between the relationship between norms on the boundary and in the elements.
\begin{lemma}\label{lemm:trace}
$\forall \bv\in V_h$, for each $K^h$, it holds that
\begin{eqnarray}
	\begin{aligned}
		\int_{\pa K^h} |\bv\cdot \bn|\leq \frac{C}{h}\int_{K^h} (|\bv| +|\nabla\cdot \bv|).
	\end{aligned}
\end{eqnarray}
\end{lemma}
\begin{proof}
Let $I=[0,1]\times [0,1]$ be a reference element and $K^h=[x_i,x_{i+1}]\times [y_j,y_{j+1}]$ is the concerned fine element. Suppose $V_I$ is the space of linear functions on $I$. For each $\bv\in V_h(K^h)$, we define a corresponding $\bv_I\in V_I$ as follows,
\begin{align*}
	\bv_I(p,q)=\bv(hp+x_i,hq+y_j).
\end{align*}
Let $D$ be a domain and two norms on $V_h(D)$ are defined below.
\begin{align*}
	\cN_{D,1}(\bv)&:=\int_{\pa D}|\bv\cdot \bn|,\quad \cN_{D,2}(\bv):=\int_{D}(|\bv|+|\nabla\cdot \bv|).
\end{align*}
Owing to the relation of norms on finite-dimensional space, we can easily obtain 
\begin{align}
	N_{I,1}(\bv_I)\leq C N_{I,2}(\bv_I).\label{ref_norm}
\end{align}
Based on the definition of $\bv_I$, one can attain that 
\begin{align}
	N_{K^h,1}(\bv)=hN_{I,1}(\bv_I),\quad N_{K^h,2}(\bv)=h^2N_{I,2}(\bv_I).\label{temp_norm}
\end{align}
Combine \eqref{ref_norm} and \eqref{temp_norm} and we can obtain the desired result.
\end{proof}
\begin{lemma}\label{lemma:stab}
Recall that $\bu_t^{h,n+1}$ is the solution to \eqref{scheme2} and we have
\begin{align}
	\|u_t^{h,n+1}\|_{\kappa_n^{-1},\Omega}&\lesssim  \|q_t\|_{L^2(\Omega)}, \label{stab_1}\\
	\|u_t^{h,n+1}\|_{L^2(\Omega)}&\lesssim  \kappa_{\text{max},n}^{\frac{1}{2}}\|q_t\|_{L^2(\Omega)}. \label{stab_2}
\end{align}
\end{lemma}
\begin{proof}
Based on the simplications made at the beginning at this subsection, \eqref{scheme6} and \eqref{scheme2} can be reduced to the following form: $\forall \bv\in V_h$ and $q\in Q_h$,
\begin{align}
	(\kappa_n^{-1}\bu_t^{h,n+1},\bv)-(p_w^{h,n+1},\nabla\cdot \bv)&=0, \label{temp5}\\
	(\nabla \cdot \bu_t^{h,n+1},q)&=(q_t,q).\label{temp6}
\end{align}
By choosing $\bv=\bu_t^{h,n+1}$ in \eqref{temp5} and $q=p_w^{h,n+1}$, we have
\begin{align}
	(\kappa_n^{-1}\bu_t^{h,n+1},\bu_t^{h,n+1})=(q_t,p_w^{h,n+1}).\label{temp7}
\end{align}
From \cite{brezzi1991variational}, the Raviart-Thomas elements satisfy the following inf-sup condition:
\begin{align*}
	\|q_h\|_{L^2(\Omega)}\lesssim \sup_{\bv_h\in V_f} \frac{\int_{\Omega}\text{div}(\bv_h)q_h}{\|v_h\|_{H(\text{div};\Omega)}},
	\quad \forall q_h\in Q_h.
\end{align*}
Using \eqref{temp5}, one can obtain 
\begin{align}
	\|p_w^{h,n+1}\|_{L^2(\Omega)}\lesssim \kappa_{\text{min},n}^{-\frac{1}{2}}\|\bu_t^{h,n+1}\|_{\kappa_n^{-1},\Omega}\lesssim \|\bu_t^{h,n+1}\|_{\kappa_n^{-1},\Omega},\label{temp8}
\end{align}
where we use $\kappa_{\text{min},n}=\lam_t(S_w^{h,n})K\leq C$ for some constant $C$. 
Applying Cauchy-Schwarz inequality in \eqref{temp7} and using \eqref{temp8}, we can obtain 
\begin{align*}
	\|u_t^{h,n+1}\|_{\kappa_n^{-1},\Omega}\lesssim \|q_t\|_{L^2(\Omega)}.	
\end{align*}
Thus, \eqref{stab_2} can be attained by applying Cauchy-Schwarz inequality once again.
\end{proof}
\begin{lemma}\label{Lemma:bdy_ele}
Based on Assumption \ref{ass:fw} and Lemma \ref{lemma:stab}, we have the following estimate on $K^h$:
\begin{eqnarray}
	\begin{aligned}
		&\left|\int_{\pa K^h}\left(f_{w}(S_{w}^{*,h,n})\bu_{t}^{h,n+1} \cdot \bn-f_{w}(S_{w}^{*,H,n})\bu_{t}^{H,n+1} \cdot \bn\right)\right|	\\
		&\leq \frac{C}{h}\int_{K^{h}}|\bu_{t}^{h,n+1}-\bu_{t}^{H,n+1}|
		+\frac{C}{h}\sum_{i=0}^{N_{\text{neigh}}}\int_{K^{h,i}}|S_w^{h,n}-S_w^{H,n}|\left( \int_{K^h}|\bu_{t}^{h,n+1}| +\left|\int_{K^h}q_t\right|\right),\label{lemma:final}
	\end{aligned}
\end{eqnarray}
where we let $K^{h,0}:=K^h$ and define $N_{neigh}$ to be the number of neighboring fine elements of $K^h$.
\end{lemma}
\begin{proof}
For each $K^h$,
\begin{eqnarray}
	\begin{aligned}
		&\left|\int_{\pa K^h}\left(f_{w}(S_{w}^{*,h,n})\bu_{t}^{h,n+1} \cdot \bn-f_{w}(S_{w}^{*,H,n})\bu_{t}^{H,n+1} \cdot \bn\right)\right|\leq\\
		&\left|\int_{\pa K^h}\left(f_{w}(S_{w}^{*,h,n})-f_{w}(S_{w}^{*,H,n})\right)\bu_{t}^{H,n+1} \cdot \bn\right|+
		\left|\int_{\pa K^h}f_{w}(S_{w}^{*,H,n})\left(\bu_{t}^{h,n+1}-\bu_{t}^{H,n+1}\right) \cdot \bn\right|. \label{temp2}
	\end{aligned}
\end{eqnarray}
To deal with the first term, we prove the following inequality.
\begin{align}
	|f(S_1)-f(S_2)|\leq C|S_1-S_2|,\label{rela:f_s}
\end{align} 
for some constant $C$. As mentioned before, we only consider the case $\lambda_{w}(S)=S^p$, where $p=1,2$. Since \eqref{rela:f_s} is trivial when $p=1$, we prove the other case as follows.
\begin{align*}
	f_w(S_1)-f_w(S_2)&=\frac{S_1^2}{S_1^2+(1-S_1)^2}-\frac{S_2^2}{S_2^2+(1-S_2)^2}
	=\frac{S_1^2(1-S_2)^2-S_2^2(1-S_1)^2}{(S_1^2+(1-S_1)^2)((S_2^2+(1-S_2)^2)}\\
	&=
	\frac{(S_1+S_2-2S_1S_2)(S_1-S_2)}{(S_1^2+(1-S_1)^2)((S_2^2+(1-S_2)^2)}.
\end{align*}
By Cauchy inequality, $S_t^2+(1-S_t)^2\geq \frac{1}{2}$, for $t=1,2$. Also, by the bounds-preserving property of saturation, we have $|S_1+S_2-2S_1S_2|\leq 4$. Thus, $f_w(S_1)-f_w(S_2)\leq 16|S_1-S_2|$, which concludes the proof of \eqref{rela:f_s}.
Choosing $S_1=S_{w}^{*,h,n}$ and $S_2=S_{w}^{*,H,n}$ and applying Assumption \ref{ass:fw}, we can obtain that on each $e\subset (\pa K^h\setminus \Gamma)$, 
\begin{align*}
	|f_{w}(S_{w}^{*,h,n})-f_{w}(S_{w}^{*,H,n})|_{e}\leq C|S_{w}^{*,h,n}-S_{w}^{*,H,n}|_e
	\leq C\left(|S_{w}^{h,n}-S_{w}^{H,n}|_{K^h}+|S_{w}^{h,n}-S_{w}^{H,n}|_{K^{h,i}}\right),
\end{align*}
where $K^h\cap K^{h,i}=e$.
For $e\subset \Gamma$, it holds that 
\begin{align*}
	|f_{w}(S_{w}^{*,h,n})-f_{w}(S_{w}^{*,H,n}|_{e}\leq C|S_{w}^{h,n}-S_{w}^{H,n}|_{K^h}.
\end{align*}
Hence, we can proceed with
\begin{align}
	\left|\int_{\pa K^h}\left(f_{w}(S_{w}^{*,h,n})-f_{w}(S_{w}^{*,H,n})\right)\bu_{t}^{h,n+1} \cdot \bn\right|
	\leq C\sum_{i=0}^{N_{\text{neigh}}}|S_w^{h,n}-S_w^{H,n}|_{K^{h,i}}\int_{\pa K^h}|\bu_{t}^{h,n+1}\cdot \bn|,\label{temp4}
\end{align}
where $K^{h,0}$ is used to denote $K^h$. From Lemma \ref{lemm:trace}, it holds that 
\begin{eqnarray*}
	\begin{aligned}
		\int_{\pa K^h} |\bu_{t}^{h,n+1}\cdot \bn|\leq \frac{C}{h}\int_{K^h} (|\bu_{t}^{h,n+1}| +|\nabla\cdot \bu_{t}^{h,n+1}|).
	\end{aligned}
\end{eqnarray*}
Choosing $q=1_{K^h}$ in \eqref{scheme6} and \eqref{pp}, we have
\begin{align*}
	\int_{K^h} \nabla \bu_{t}^{h,n+1}=	\int_{K^h} \nabla \bu_{t}^{H,n+1}=\int_{K^h} q_t.
\end{align*}
Since both $\nabla \cdot \bu_{t}^{h,n+1}$ and $\nabla \cdot \bu_{t}^{H,n+1}$  are constant in $K^h$, we can attain 
\begin{align}
	|K^h|\nabla \cdot\bu_{t}^{h,n+1}=|K^h|\nabla \cdot\bu_{t}^{H,n+1}=\int_{K^h} q_t,\label{nabla_u}
\end{align}
where $|K^h|$ is the area of $K^h$.
Hence, \eqref{temp4} can be replaced by 
\begin{align}
	\left|\int_{\pa K^h}\left(f_{w}(S_{w}^{*,h,n})-f_{w}(S_{w}^{*,H,n})\right)\bu_{t}^{h,n+1} \cdot \bn\right|
	\leq \frac{C}{h}\sum_{i=0}^{N_{\text{neigh}}}|S_w^{h,n}-S_w^{H,n}|_{K^{h,i}}\left( \int_{K^h}|\bu_{t}^{h,n+1}| +\left|\int_{K^h}q_t\right|\right).\label{temp4_2}
\end{align}
Then we can deal with the second term on the right hand side of \eqref{temp2}. Using the relation proved in Lemma \ref{lemm:trace}, 
\begin{eqnarray}
	\begin{aligned}
		&\left|\int_{\pa K^h}f_{w}(S_{w}^{*,H,n})\left(\bu_{t}^{h,n+1}-\bu_{t}^{H,n+1}\right) \cdot \bn\right|\leq
		\int_{\pa K^h}\left|\left(\bu_{t}^{h,n+1}-\bu_{t}^{H,n+1}\right) \cdot \bn\right|\\
		&\leq \frac{C}{h}\int_{K^h}\left(|\bu_{t}^{h,n+1}-\bu_{t}^{H,n+1}|+|\nabla\cdot(\bu_{t}^{h,n+1}-\bu_{t}^{H,n+1})|\right).
		\label{temp3}
	\end{aligned}
\end{eqnarray}
Applying \eqref{nabla_u} again, we can transform \eqref{temp3} to
\begin{align}
	\left|\int_{\pa K^h}f_{w}(S_{w}^{*,H,n})\left(\bu_{t}^{h,n+1}-\bu_{t}^{H,n+1}\right) \cdot \bn\right|\leq 
	\frac{C}{h}\int_{K^h}|\bu_{t}^{h,n+1}-\bu_{t}^{H,n+1}|.\label{temp3_2}
\end{align}
Finally, by combining \eqref{temp2}, \eqref{temp4_2} and \eqref{temp3_2}, we conclude with \eqref{lemma:final}.
\end{proof}
In the following two lemmas, we give estimations on  $\|\bu_t^{h,n+1}-\bu_t^{H,n+1}\|_{L^2(\Omega)}$.
Before this, we present some notations.
\begin{enumerate}
\item
$C_{\text{err}}=\frac{C_VH}{h}$ and $C_V$ depends on the polynomial order of the fine-grid bases in $V_{\text{snap}}$.
\item 
$R_{D}^{(i)}$ is the residual operator after $i$ iterations at time $t_{p}$, which is defined in \eqref{loc_res}. $B_i$ is the number of non-overlapping domains, on which we perform the bases enrichment.
\end{enumerate}

Then, we estimate the velocity error $\|\bu_t^{h,n+1}-\bu_t^{H,n+1}\|_{L^2(\Omega)}$. We split the difference in four parts. In particular, we define 
\begin{align*}
\bw_1&=(\kappa_n^H)^{-\frac{1}{2}}\bu_t^{h,n+1}-\kappa_n^{-\frac{1}{2}}\bu_t^{h,n+1},\\
\bw_2&=\kappa_n^{-\frac{1}{2}}\bu_t^{h,n+1}-\kappa_0^{-\frac{1}{2}}\bu_t^{h,1},\\
\bw_3&=\kappa_0^{-\frac{1}{2}}\bu_t^{h,1}-\kappa_0^{-\frac{1}{2}}\bu_t^{H,1},\\
\bw_4&=\kappa_0^{-\frac{1}{2}}\bu_t^{H,1}-(\kappa_n^H)^{-\frac{1}{2}}\bu_t^{H,n+1}.
\end{align*}
Then 
\begin{align}
\bu_t^{h,n+1}-\bu_t^{H,n+1}=(\kappa_n^H)^{\frac{1}{2}}\sum_{i=1}^4\bw_i.\label{err_split}
\end{align} 
We remark that multiscale bases constructed at initial time are used in different time steps. Meanwhile, $\kappa_n$ and $\kappa_n^H$ are dependent on time $t_n$, which are different with $\kappa_{0}$ at the beginning. As a result, $\bw_2$ and $\bw_4$ should be considered. As for  $\bw_1$ and $\bw_4$, they are resulted from approximation with multiscale space. In the following two lemmas, we focus on estimate $\|\bw_2\|_{L^2(\Omega)}$ to $\|\bw_4\|_{L^2(\Omega)}$.
\begin{lemma}\label{lemma:e_2}
Recall that  $\bu_t^{h,n+1}$ and $\bu_t^{h,1}$ are solutions to \eqref{scheme2} at time $t_{n+1}$ and $t_1$. Correspondingly, $\bu_t^{H,n+1}$ and $\bu_t^{H,1}$ are solutions to \eqref{ms-scheme3} at time $t_{n+1}$ and $t_1$. Then the following estimations hold,
\begin{align}
	\|\kappa_n^{-\frac{1}{2}}\bu_t^{h,n+1}-\kappa_0^{-\frac{1}{2}}\bu_t^{h,1}\|_{L^2(\Omega)}&\leq C \|S_w^{h,n}-S_w^{h,0}\|_{L^{\infty}(\Omega)}\|q_t\|_{L^2(\Omega)}, \label{e_2}\\
	\|\kappa_0^{-\frac{1}{2}}\bu_t^{H,1}-(\kappa_n^H)^{-\frac{1}{2}}\bu_t^{H,n+1}\|_{L^2(\Omega)}&\leq C \|S_w^{H,n}-S_w^{H,0}\|_{L^{\infty}(\Omega)}\|q_t\|_{L^2(\Omega)}, \label{e_4}
\end{align}
\end{lemma} 
\begin{proof}
We first recall  \eqref{scheme2} and \eqref{scheme6} at time $t_{n+1}$ and $t_1$.
\begin{align}
	\int_\Omega  \kappa_{n}^{-1}\bu_t^{h,n+1}\cdot\bv-\int_\Omega\nabla\cdot\bv p_w^{h,n+1}&=0,\quad\forall \bv\in V_f(\Omega),\label{1_1} \\ 
	\int_\Omega  \nabla\cdot\bu_t^{h,n+1}q&=\int_\Omega  q_tq,\quad\forall q\in Q_f.\label{1_2} 
\end{align}
\begin{align}
	\int_\Omega  \kappa_0^{-1}\bu_t^{h,1}\cdot\bv-\int_\Omega\nabla\cdot\bv  p_w^{h,1}&=0,\quad\forall \bv\in V_f(\Omega), \\ 
	\int_\Omega  \nabla\cdot\bu_t^{h,n+1}q&=\int_\Omega  q_tq,\quad\forall q\in Q_f.\label{1_3}
\end{align}
Taking $\bv=\bu_t^{h,n+1}$ in \eqref{1_1} and $q=p_h^{h,n+1}$ in \eqref{1_2}, it holds that
\begin{align}
	\int_\Omega  \kappa_n^{-1}|\bu_t^{h,n+1}|^2-\int_\Omega\nabla \cdot\bu_t^{h,n+1}p_w^{h,n+1}&=0, \\ 
	\int_\Omega  \text{div}(\bu_t^{h,n+1})p_w^{h,n+1}&=\int_\Omega  q_tp_w^{h,n+1}.
\end{align}
It follows that
\begin{align}
	\int_\Omega  \kappa_n^{-1}|\bu_t^{h,n+1}|^2=\int_\Omega  q_tp_w^{h,n+1}. 
\end{align}
Choosing $\bv=\bu_t^{h,1}$ in \eqref{1_1} and $q=p_w^{h,n+1}$ in \eqref{1_3}, we obtain
\begin{align}
	\int_\Omega  \kappa_n^{-1}\bu_t^{h,n+1}\cdot\bu_t^{h,1}=\int_\Omega  q_tp_w^{h,n+1}. 
\end{align}
Thus,
\begin{align}
	\int_\Omega  \kappa_n^{-1}\bu_t^{h,n+1}\cdot\bu_t^{h,1}=\int_\Omega  \kappa_n^{-1}|\bu_t^{h,n+1}|^2.\label{1_4} 
\end{align}
Similarly, one can prove that 
\begin{align}
	\int_\Omega  \kappa_0^{-1}\bu_t^{h,1}\cdot\bu_t^{h,n+1}=\int_\Omega  \kappa_0^{-1}|\bu_t^{h,1}|^2. \label{1_5} 
\end{align}
Using \eqref{1_4} and \eqref{1_5}, it holds that
\begin{align*}
	\int_{\Omega}|\kappa_0^{-\frac{1}{2}}\bu_t^{h,1}-\kappa_n^{-\frac{1}{2}}\bu_t^{h,n+1}|^2&=\int_{\Omega}\left(|\kappa_0^{-\frac{1}{2}}\bu_t^{h,1}|^2+|\kappa_n^{-\frac{1}{2}}\bu_t^{h,n+1}|^2-2(\kappa_0\kappa_n)^{-\frac{1}{2}}\right)\bu_t^{h,1}\cdot \bu_t^{h,n+1}\\
	&=\int_{\Omega}|\kappa_0^{-\frac{1}{2}}-\kappa_n^{-\frac{1}{2}}|^2\bu_t^{h,1}\cdot \bu_t^{h,n+1}\\
	&=\int_{\Omega}\left|\left(\frac{\kappa_{0}}{\kappa_n}\right)^{\frac{1}{4}}-\left(\frac{\kappa_{n}}{\kappa_0}\right)^{\frac{1}{4}}\right|^2(\kappa_{0}^{-\frac{1}{2}}\bu_t^{h,1})\cdot (\kappa_{n}^{-\frac{1}{2}}\bu_t^{h,n+1}).
\end{align*}
Since
\begin{align*}
	\left|\left(\frac{\kappa_{0}}{\kappa_n}\right)^{\frac{1}{4}}-\left(\frac{\kappa_{n}}{\kappa_0}\right)^{\frac{1}{4}}\right|\leq C\left|\left(\frac{\kappa_{0}}{\kappa_n}\right)-\left(\frac{\kappa_{n}}{\kappa_0}\right)\right|\leq C|S_w^{h,0}-S_w^{h,n}|,
\end{align*}
we conclude that 
\begin{align*}
	\|\kappa_0^{-\frac{1}{2}}\bu_t^{h,1}-\kappa_n^{-\frac{1}{2}}\bu_t^{h,n+1}\|_{L^2(\Omega)}^2&\leq C\|S_w^{h,0}-S_w^{h,n}\|_{L^{\infty}(\Omega)}^2\|\bu_t^{h,1}\|_{\kappa_{0}^{-1},\Omega}\|\bu_t^{h,n+1}\|_{\kappa_{0}^{-1},\Omega},\\
	&\leq C\|S_w^{h,0}-S_w^{h,n}\|_{L^{\infty}(\Omega)}^2\|q_t\|_{L^2(\Omega)}^2,
\end{align*}
where the last inequality holds based on Lemma \ref{lemma:stab}. Similarly, \eqref{e_4} can be obtained.
\end{proof} 
Next we evaluate $\|\bw_3\|_{L^2(\Omega)}$.
\begin{lemma} \label{lemma:e_3}
\cite{online-mixed}We recall that $\bu_{t}^{h,1}$ and $\bu_t^{H,1}$ are solutions to \eqref{scheme2} and  \eqref{ms-scheme3} at time $t_1$, respectively. Then we have the following estimation:
\begin{align*}
	\|\bu_{t}^{h,1}-\bu_{t}^{H,1}\|_{\kappa_0^{-1},\Omega}\leq C_{\text{err}}\sum_{j=1}^{N_{\text{in},c}}\|R_{D_j}^{(0)}\|_{(V_{\text{snap}}^{j})^{*}}(\lambda_j^{l_j+1})^{-1}-
	\sum_{i=1}^{m}\sum_{j=1}^{N_{\text{in},c}}\|R_{D_j}^{(i)}\|_{(\hat{V}_{D_j})^{*}},
\end{align*}
where $C_{\text{err}}=\frac{C_V H}{h}$ and $C_V$ depends on the polynomial order of the fine-grid basis functions in $V_{\text{snap}}$.
\end{lemma}
Based on above preparations, we conclude the error estimation for velocity in the theorem below.
\begin{theorem}\label{lemma:velo_error}
We recall that $\bu_t^{h,n+1}$ and $\bu_t^{H,n+1}$ are solutions to \eqref{scheme2} and \eqref{ms-scheme3} at time $t_{n+1}$. Suppose we have used $m$ iterations of bases enrichment at time $t_1$. Then  \begin{eqnarray}
	\begin{aligned}
		\|\bu_t^{h,n+1}-\bu_t^{H,n+1}\|_{L^2(\Omega)}&\leq
		C\|S_w^{h,n}-S_w^{H,n}\|_{L^2(\Omega)}+
		C(\kappa_{\text{max},n}^H)^{\frac{1}{2}}\left(C_{\text{err}}\sum_{j=1}^{N_{\text{in},c}}\|R_{D_j}^{(0)}\|_{(V_{\text{snap}}^{j})^{*}}(\lambda_{l_j+1}^{(j)})^{-1}\right.\\
		&\left.+\|S_w^{h,n}-S_w^{h,0}\|_{L^{\infty}(\Omega)}+\|S_w^{H,n}-S_w^{H,0}\|_{L^{\infty}(\Omega)}-
		\sum_{i=1}^{m}\sum_{j=1}^{N_{\text{in},c}}\|R_{D_j}^{(i)}\|_{(\hat{V}_{D_j})^{*}}\right).\label{lem:result_velo}
	\end{aligned}
\end{eqnarray}
\end{theorem}
\begin{proof}
The goal of this part is to evaluate $\|\bu_t^{h,n+1}-\bu_t^{H,n+1}\|_{L^2(\Omega)}$. To this end, we use an error splitting as follows.
\begin{eqnarray}
	\begin{aligned}
		\bu_t^{h,n+1}-\bu_t^{H,n+1}&=(\kappa_n^H)^{\frac{1}{2}}\left((\kappa_n^H)^{-\frac{1}{2}}\bu_t^{h,n+1}-(\kappa_n^H)^{-\frac{1}{2}}\bu_t^{H,n+1}\right)\\
		&=(\kappa_n^H)^{\frac{1}{2}}\left((\kappa_n^H)^{-\frac{1}{2}}\bu_t^{h,n+1}-\kappa_n^{-\frac{1}{2}}\bu_t^{h,n+1}+\kappa_n^{-\frac{1}{2}}\bu_t^{h,n+1}-(\kappa_n^H)^{-\frac{1}{2}}\bu_t^{H,n+1}\right).\label{v_diff1}
	\end{aligned}	
\end{eqnarray}
Then, we deal with two parts of the right hand side in \eqref{v_diff1} separately.
\begin{eqnarray}
	\begin{aligned}
		\left\|(\kappa_n^H)^{\frac{1}{2}}\left((\kappa_n^H)^{-\frac{1}{2}}\bu_t^{h,n+1}-\kappa_n^{-\frac{1}{2}}\bu_t^{h,n+1}\right)\right\|_{L^2(\Omega)}\leq
		\kappa_{\text{max},n}^{\frac{1}{2}}\left\|1-\left(\frac{\kappa_n^H}{\kappa_n}\right)^{\frac{1}{2}}\right\|_{L^2(\Omega)}\|\bu_t^{h,n+1}\|_{\kappa_n^{-1},\Omega}.\label{v_diff2}
	\end{aligned}	
\end{eqnarray}
By the definitions of $\kappa_n$ and $\kappa_n^H$,  we have
\begin{align*}
	\left|1-\left(\frac{\kappa_n^H}{\kappa_n}\right)^{\frac{1}{2}}\right|&=\left|1-\left(\frac{\lambda_t(S_w^{H,n})}{\lambda_t(S_w^{h,n})}\right)^{\frac{1}{2}}\right|
	=\left|\frac{(\lambda_t(S_w^{h,n}))^{\frac{1}{2}}-(\lambda_t(S_w^{H,n}))^{\frac{1}{2}}}{(\lambda_t(S_w^{h,n}))^{\frac{1}{2}}}\right|\\
	&=\left|\frac{\lambda_t(S_w^{h,n})-\lambda_t(S_w^{H,n})}{(\lambda_t(S_w^{h,n}))^{\frac{1}{2}}\left((\lambda_t(S_w^{h,n}))^{\frac{1}{2}}+(\lambda_t(S_w^{H,n}))^{\frac{1}{2}}\right)}\right|\\
	&\leq C|S_w^{h,n}-S_w^{H,n}|.
\end{align*} 
From Lemma 
\ref{lemma:stab}, \eqref{v_diff2} can be replaced by
\begin{align}
	\left\|(\kappa_n^H)^{\frac{1}{2}}\left((\kappa_n^H)^{-\frac{1}{2}}\bu_t^{h,n+1}-\kappa_n^{-\frac{1}{2}}\bu_t^{h,n+1}\right)\right\|_{L^2(\Omega)}\leq	C\kappa_{\text{max},n}^{\frac{1}{2}}\|S_w^{h,n}-S_w^{H,n}\|_{L^2(\Omega)}.\label{v_diff3}
\end{align}
For the other part in \eqref{v_diff1}, we have
\begin{align*}
	&\|\kappa_n^{-\frac{1}{2}}\bu_t^{h,n+1}-(\kappa_n^H)^{-\frac{1}{2}}\bu_t^{H,n+1}\|_{L^2(\Omega)}\leq \\
	&\|\kappa_n^{-\frac{1}{2}}\bu_t^{h,n+1}-\kappa_0^{-\frac{1}{2}}\bu_t^{h,1}\|_{L^2(\Omega)}+\|\kappa_0^{-\frac{1}{2}}\bu_t^{h,1}-\kappa_0^{-\frac{1}{2}}\bu_t^{H,1}\|_{L^2(\Omega)}+\|\kappa_0^{-\frac{1}{2}}\bu_t^{H,1}-(\kappa_n^H)^{-\frac{1}{2}}\bu_t^{H,n+1}\|_{L^2(\Omega)}\\
	&:=e_1+e_2+e_3.
\end{align*}
Summarize the conclusion in Lemma \ref{lemma:e_2} to \ref{lemma:e_3} and we have
\begin{align*}
	&e_1\leq C \|S_w^{h,n}-S_w^{h,0}\|_{L^{\infty}(\Omega)}\|q_t\|_{L^2(\Omega)},\\
	&e_2\leq C_{\text{err}}\sum_{j=1}^{N_{\text{in},c}}\|R_{D_j}^{(0)}\|_{(V_{\text{snap}}^{j})^{*}}(\lambda_{l_j+1}^{(j)})^{-1}-
	\sum_{i=1}^{m}\sum_{j=1}^{N_{\text{in},c}}\|R_{D_j}^{(i)}\|_{(\hat{V}_{D_j})^{*}},\\
	&e_3\leq C \|S_w^{H,n}-S_w^{H,0}\|_{L^{\infty}(\Omega)}\|q_t\|_{L^2(\Omega)},
\end{align*}
where $m$ is the total number of iterations of bases enrichment.
Thus, it arrives at 
\begin{eqnarray}
	\begin{aligned}
		\|\kappa_n^{-\frac{1}{2}}\bu_t^{h,n+1}-(\kappa_n^H)^{-\frac{1}{2}}\bu_t^{H,n+1}\|_{L^2(\Omega)}&\leq
		C\left(\|S_w^{h,n}-S_w^{h,0}\|_{L^{\infty}(\Omega)}+\|S_w^{H,n}-S_w^{H,0}\|_{L^{\infty}(\Omega)}\right)\\
		&+C_{\text{err}}\sum_{j=1}^{N_{\text{in},c}}\|R_{D_j}^{(0)}\|_{(V_{\text{snap}}^{j})^{*}}(\lambda_{l_j+1}^{(j)})^{-1}-
		\sum_{i=1}^{m}\sum_{j=1}^{N_{\text{in},c}}\|R_{D_j}^{(i)}\|_{(\hat{V}_{D_j})^{*}}.\label{v_diff4}
	\end{aligned}
\end{eqnarray}
Combining \eqref{v_diff3} and \eqref{v_diff4}, we conclude with \eqref{lem:result_velo}.
\end{proof}
Finally, we estimate the saturation error.
\begin{theorem}
Recall that $S_w^{h,n+1}$ and $S_w^{H,n+1}$ are solutions to \eqref{scheme1} and \eqref{ms-scheme2}, respectively. We have the following estimation.
\begin{eqnarray}
	\begin{aligned}
		\|S_{w}^{h,n+1}-S_{w}^{H,n+1}\|_{L^2(\Omega)}^2\leq  C_1(\Delta t)
		\sum_{i=0}^{n}(C_2(\Delta t))^{n-i}(e_{\bv}^{i+1})^2,
		\label{final_result}
	\end{aligned}
\end{eqnarray}
where  $C_1(\Delta t)=\frac{C\Delta t}{h(\zeta-C\Delta t)}$, $C_2(\Delta t)=(\frac{h\zeta}{\Delta t}+C_{\kappa}+1 )C_1(\Delta t)$ , and 
\begin{align*}
	C_{\kappa}:=&4(1+\kappa_{\text{max},0})\|q_t\|_{L^2(\Omega)}^2,\\
	e_{\bv}^{i+1}:=&
	C(\kappa_{\text{max},i}^H)^{\frac{1}{2}}\left(C_{\text{err}}\sum_{j=1}^{N_{\text{in},c}}\|R_{D_j}^{(0)}\|_{(V_{\text{snap}}^{j})^{*}}(\lambda_{l_j+1}^{(j)})^{-1}+\|S_w^{h,i}-S_w^{h,0}\|_{L^{\infty}(\Omega)}\right.\\
	&\left.+\|S_w^{H,i}-S_w^{H,0}\|_{L^{\infty}(\Omega)}-
	\sum_{p=1}^{m}\sum_{j=1}^{N_{\text{in},c}}\|R_{D_j}^{(p)}\|_{(\hat{V}_{D_j})^{*}}\right).
\end{align*}
\end{theorem}
\begin{proof}
First we rewrite the fine-scale and coarse-scale transport equations. $\forall q\in Q_h$,
\begin{align}
	&(\zeta\dfrac{S_{w}^{h,n+1}-S_{w}^{h,n}}{t_{n+1}-t_n},q)+\beta_{w}(\bu_t^{h,n+1},q;S_w^{h,n})=(q_{w},q),\label{ref-Scheme1}	\\
	&(\zeta\dfrac{S_{w}^{H,n+1}-S_{w}^{H,n}}{t_{n+1}-t_n},q)+\beta_{w}(\bu_t^{H,n+1},q;S_w^{H,n})=(q_{w},q),\label{ms-Scheme1}
\end{align}
where for $\alpha=w,n$, $\beta_{\alpha}$ is defined in \eqref{beta_alpha}.
Subtract \eqref{ref-Scheme1} by \eqref{ms-Scheme1} and we attain 
\begin{equation}
	(S_{w}^{h,n+1}-S_{w}^{H,n+1},q)-(S_{w}^{h,n}-S_{w}^{H,n},q)=-\frac{\Delta t}{\zeta}\left(\beta_{w}(\bu_t^{h,n+1},q;S_w^{h,n})-\beta_{w}(\bu_t^{H,n+1},q;S_w^{H,n})\right). \label{differ}
\end{equation}
For each $K^h\in \cT_h$, we can estimate $\|S_{w}^{h,n+1}-S_{w}^{H,n+1}\|_{L^2(K^h)}^2$ separately. To this end, we choose $q=(S_{w}^{h,n+1}-S_{w}^{H,n+1})1_{K^h}$. 
Hence \eqref{differ} can be transformed to 
\begin{eqnarray}
	\begin{aligned}
		&\|S_{w}^{h,n+1}-S_{w}^{H,n+1}\|_{L^2(K^h)}^2-\int_{K^h}(S_{w}^{h,n+1}-S_{w}^{H,n+1})(S_{w}^{h,n}-S_{w}^{H,n})\\
		&=-\frac{\Delta t}{\zeta} (S_{w}^{h,n+1}-S_{w}^{H,n+1})_{K^h}\int_{\pa K^h}\left(f_{w}(S_{w}^{*,h,n})\bu_{t}^{h,n+1} \cdot \bn-f_{w}(S_{w}^{*,H,n})\bu_{t}^{H,n+1} \cdot \bn\right),\label{differ2}
	\end{aligned}
\end{eqnarray}
Apply the result in Lemma \ref{Lemma:bdy_ele} and we can attain that 
\begin{eqnarray}
	\begin{aligned}
		&\|S_{w}^{h,n+1}-S_{w}^{H,n+1}\|_{L^2(K^h)}^2-\int_{K^h}(S_{w}^{h,n+1}-S_{w}^{H,n+1})(S_{w}^{h,n}-S_{w}^{H,n})\\
		&\leq C\frac{\Delta t}{\zeta} (S_{w}^{h,n+1}-S_{w}^{H,n+1})_{K^h} \left(\frac{1}{h}\int_{K^h}|\bu_t^{h,n+1}-\bu_t^{H,n+1}|^2+\frac{1}{h}\sum_{i=0}^{N_{\text{neigh}}}\int_{K^{h,i}}|S_w^{h,n}-S_w^{H,n}|\left( \int_{K^h}|\bu_{t}^{h,n+1}| +\left|\int_{K^h}q_t\right|\right)\right)\\
		&\leq \frac{C\Delta t}{2\zeta}\left(6\|S_{w}^{h,n+1}-S_{w}^{H,n+1}\|_{L^2(K^h)}^2+\frac{1}{h}\sum_{i=0}^{N_{\text{neigh}}}\|S_{w}^{h,n}-S_{w}^{H,n}\|_{L^2(K^{h,i})}^2\left( \|\bu_{t}^{h,n+1}\|_{L^2(K^{h,i})}^2 +\|q_t\|_{L^2(K^{h,i})}^2\right)\right.\\
		&\left.+\frac{1}{h}\|\bu_t^{h,n+1}-\bu_t^{H,n+1}\|_{L^2(K^{h})}^2\right),
		\label{differ3}
	\end{aligned}
\end{eqnarray}
where we use Cauchy-Schwarz inequality in the last step. Here we recall that $K^{h,0}:=K^h$ and define $N_{neigh}$ to be the number of neighboring fine elements of $K^h$.
Summing over all $K^h\in \Omega$ and applying Cauchy-Schwarz inequality again , it holds that 
\begin{eqnarray}
	\begin{aligned}
		&\|S_{w}^{h,n+1}-S_{w}^{H,n+1}\|_{L^2(\Omega)}^2-\int_{\Omega}(S_{w}^{h,n+1}-S_{w}^{H,n+1})(S_{w}^{h,n}-S_{w}^{H,n})\\
		&\leq \frac{C\Delta t}{2\zeta}\left(\|S_{w}^{h,n+1}-S_{w}^{H,n+1}\|_{L^2(\Omega)}^2+\frac{1}{h}\|S_{w}^{h,n}-S_{w}^{H,n}\|_{L^2(\Omega)}^2\left( \|\bu_{t}^{h,n+1}\|_{L^2(\Omega)}^2 +\|q_t\|_{L^2(\Omega)}^2\right)+\frac{1}{h}\|\bu_t^{h,n+1}-\bu_t^{H,n+1}\|_{L^2(\Omega)}^2\right)\\
		&\leq \frac{C\Delta t}{2\zeta}\left(\|S_{w}^{h,n+1}-S_{w}^{H,n+1}\|_{L^2(\Omega)}^2+\frac{C_{\kappa}}{h}\|S_{w}^{h,n}-S_{w}^{H,n}\|_{L^2(\Omega)}^2+\frac{1}{h}\|\bu_t^{h,n+1}-\bu_t^{H,n+1}\|_{L^2(\Omega)}^2\right),
		\label{differ4}
	\end{aligned}
\end{eqnarray}
where $C_{\kappa}:=8(1+\kappa_{\text{max},0})\|q_t\|_{L^2(\Omega)}^2$. Here we use $\kappa_{\text{max},n}\leq 2\kappa_{\text{max},0}$ for all $n$.
Applying Cauchy-Schwarz inequality once again, we obtain that
\begin{eqnarray}
	\begin{aligned}
		\|S_{w}^{h,n+1}-S_{w}^{H,n+1}\|_{L^2(\Omega)}^2\leq
		\frac{\zeta h+CC_{\kappa}\Delta t}{h(\zeta-C\Delta t)}\|S_{w}^{h,n}-S_{w}^{H,n}\|_{L^2(\Omega)}^2+\frac{C\Delta t }{h(\zeta-C\Delta t)}\|\bu_t^{h,n+1}-\bu_t^{H,n+1}\|_{L^2(\Omega)}^2.
		\label{differ5}
	\end{aligned}
\end{eqnarray}
From Theorem \ref{lemma:velo_error}, we have
\begin{eqnarray*}
	\begin{aligned}
		\|\bu_t^{h,n+1}-\bu_t^{H,n+1}\|_{L^2(\Omega)}&\leq
		C\|S_w^{h,n}-S_w^{H,n}\|_{L^2(\Omega)}+
		C(\kappa_{\text{max},n}^H)^{\frac{1}{2}}\left(C_{\text{err}}\sum_{j=1}^{N_{\text{in},c}}\|R_{D_j}^{(0)}\|_{(V_{\text{snap}}^{j})^{*}}(\lambda_{l_j+1}^{(j)})^{-1}\right.\\
		&\left.+\|S_w^{h,n}-S_w^{h,0}\|_{L^{\infty}(\Omega)}+\|S_w^{H,n}-S_w^{H,0}\|_{L^{\infty}(\Omega)}-
		\sum_{i=1}^{m}\sum_{j=1}^{N_{\text{in},c}}\|R_{D_j}^{(i)}\|_{(\hat{V}_{D_j})^{*}}\right).
	\end{aligned}
\end{eqnarray*}
Then, we define  
\begin{align*}
	e_{\bv}^{n+1}&:=
	C(\kappa_{\text{max},n}^H)^{\frac{1}{2}}\left(C_{\text{err}}\sum_{j=1}^{N_{\text{in},c}}\|R_{D_j}^{(0)}\|_{(V_{\text{snap}}^{j})^{*}}(\lambda_{l_j+1}^{(j)})^{-1}+
	\|S_w^{h,n}-S_w^{h,0}\|_{L^{\infty}}\right.\\
	&\left.+\|S_w^{H,n}-S_w^{H,0}\|_{L^{\infty}} -
	\sum_{i=1}^{m}\sum_{j=1}^{N_{\text{in},c}}\|R_{D_j}^{(i)}\|_{(\hat{V}_{D_j})^{*}}\right)
\end{align*}
Then
\begin{eqnarray}
	\begin{aligned}
		\|S_{w}^{h,n+1}-S_{w}^{H,n+1}\|_{L^2(\Omega)}^2\leq
		\left(\frac{\zeta h+CC_{\kappa}\Delta t}{h(\zeta-C\Delta t)}+\frac{C\Delta t }{h(\zeta-C\Delta t)}\right)\|S_{w}^{h,n}-S_{w}^{H,n}\|_{L^2(\Omega)}^2+\frac{C\Delta t }{h(\zeta-C\Delta t)}(e_{\bv}^{n+1})^2.
		\label{differ5-2}
	\end{aligned}
\end{eqnarray}
Let $C_1(\Delta t)=\frac{C\Delta t}{h(\zeta-C\Delta t)}$ and $C_2(\Delta t)=(\frac{h\zeta}{\Delta t}+C_{\kappa}+1 )C_1(\Delta t)$ , then \eqref{differ5-2} can be reduced to 
\begin{eqnarray}
	\begin{aligned}
		\|S_{w}^{h,n+1}-S_{w}^{H,n+1}\|_{L^2(\Omega)}^2\leq
		C_2(\Delta t)\|S_{w}^{h,n}-S_{w}^{H,n}\|_{L^2(\Omega)}^2+ C_1(\Delta t)(e_{\bv}^{n+1})^2.
		\label{differ6}
	\end{aligned}
\end{eqnarray}
Define $P_h$ be the $L^2-$ projection into $Q_h$ and based on $S_{w}^{h,0}=S_{w}^{h,0}=P_h(S_w^0)$,
we finish the error estimation for saturation with
\begin{eqnarray}
	\begin{aligned}
		\|S_{w}^{h,n+1}-S_{w}^{H,n+1}\|_{L^2(\Omega)}^2\leq  C_1(\Delta t)
		\sum_{i=1}^{n+1}(C_2(\Delta t))^{n+1-i}(e_{\bv}^{i})^2.
		\label{differ7}
	\end{aligned}
\end{eqnarray}
Indeed, since \eqref{differ6} is relatively simple, induction is sufficient to prove \eqref{differ7} and it is unnecessary to apply Discrete Gronwall inequality \cite{holte2009discrete}. 
\end{proof}

\section{Conclusion}
In this work, we propose a physics-preserving multiscale scheme to solve a two-phase flow problem with heterogeneous media. The scheme can satisfy the local conservation for two phases and it is proved to be unbiased and conditional bounds-preserving for saturation, where  a sufficiently small time step is required. Moreover, the scheme is well-designed because of using multiscale bases. In particular,  one solve the saturation on a fine mesh while the velocity is computed on a coarse mesh by utilizing a set of multiscale bases. It is notable that the multiscale bases are computed on a fine mesh to capture the fine-scale information contained in the permeability field. Besides, we use some residual-driven  bases to further reduce the velocity error and those new bases are computed on a fine mesh as well. Hence, contributed by the representative approximation space, this method can achieve a good computation efficiency but not lose much precision.

To illustrate the performance of the method, we use a set of well-designed numerical experiments. Through these experiments, one can see the proposed scheme is physics-preserving, i.e., the local conservation of mass for both phases are satisfied. Moreover, the error is convergent with smaller time step and coarse mesh size. Besides, the effect of using more bases is also distinguishable especially when residual-driven bases are constructed. Hence, from  numerical results, the method is efficient and relatively accurate. On the other hand, we provide a rigorous analysis which serves as the base of the numerical performance. The analysis roughly contains two parts. We first prove the method is physics-preserving, i.e., it is locally conservative for both phases. Besides, we prove that it is unbiased as well as conditional bounds-preserving, where the condition is to choose a small enough time step. In the second part, we prove the error convergence, which is based on exploring the relation of velocity error and saturation error. From the final conclusion, one can observe that the saturation  error indeed depends on time step size, coarse mesh size as well as the number of multiscale bases. In addition, shown in the final conclusion, the error convergence of multiscale bases in two stages has different forms, which sheds some
light on the fact that residual-driven bases are more powerful in reducing errors. Overall, the numerical results are consistent to analysis. From both perspectives of analysis and numerical examples, the proposed method is efficient as well as relevant accurate.

\section*{Acknowledgement}

The research of Eric Chung is partially supported by the Hong Kong RGC General Research Fund (Project: 14304021). 
The research of Shuyu Sun is partially supported by King Abdullah University of Science and Technology (KAUST) through the grants BAS/1/1351-01, URF/1/4074-01, and URF/1/3769-01.
	\bibliographystyle{plain}
	\bibliography{references}
\end{document}